\documentclass[11pt]{amsart}
\allowdisplaybreaks[1]

\usepackage{a4wide}
\usepackage{amssymb}
\usepackage{amsmath}

\newtheorem{theorem}{Theorem}
\newtheorem{lemma}{Lemma}
\newtheorem{prop}{Proposition}
\newtheorem{coro}{Corollary}
\newtheorem{definition}{Definition}

\newcommand{\RR}{\mathbb{R}}

\newcommand{\CC}{\mathbb{C}}
\newcommand{\NN}{\mathbb{N}}

\newcommand{\dd}{\,{\rm d}}

\newcommand{\UKV}{U_{K,V}}

\newcommand{\MM}{\mathcal{M}(G)}
\newcommand{\MCV}{\mathcal{M}_{C,V}(G)}
\newcommand{\MTB}{\mathcal{M}^{\infty}(G)}

\newcommand{\oO}{(\varOmega,\alpha)}
\newcommand{\cHp}{\mathcal{H}_{pp}(T)}

\newcommand{\cU}{\hspace{1pt}\mathcal{U}}
\newcommand{\cV}{\mathcal{V}}

\newcommand{\cUD}{\hspace{1pt}\mathcal{U D}(G)}
\newcommand{\cDV}{\mathcal{D}_V (G)}
\newcommand{\cC}{\mathcal{C}(G)}
\newcommand{\cS}{\mathcal{S}}

\newcommand{\UFell}{\mathcal{U}(C,\mathcal{F})}
\newcommand{\Uwir}{U_{K,V} (H)}  
\newcommand{\CalF}{{\mathcal{F}}}

\DeclareMathOperator{\supp}{supp}

\begin{document}

\title[Pure point spectrum]
{Dynamical systems on translation bounded measures: \\[1mm]
Pure point dynamical and diffraction spectra}

\author{Michael Baake}
\address{Fakult\"at f\"ur Mathematik, Universit\"at
Bielefeld, Postfach 100131, 33501 Bielefeld, Germany}
\email{mbaake@mathematik.uni-bielefeld.de}

\author{Daniel Lenz}
\address{Fakult\"at  f\"ur Mathematik, TU Chemnitz,
09107 Chemnitz, Germany}
\email{dlenz@mathematik.tu-chemnitz.de}

\begin{abstract}  
Certain topological dynamical systems are considered
that arise from actions of $\sigma$-compact locally compact Abelian 
groups on compact spaces of translation bounded measures. 
Such a measure dynamical system is shown to have pure point 
dynamical spectrum if and only if its diffraction spectrum is pure point.
\end{abstract}
\maketitle

\section{Introduction}
This paper deals with certain dynamical systems build from measures on
$\sigma$-compact locally compact Abelian groups.  These dynamical
systems give rise to two spectra: the dynamical spectrum and the
diffraction spectrum. After introducing the dynamical systems and
discussing their basic topological features, we will focus on studying
the relationship between these two spectra. Particular attention will
be paid to the case where one of the spectra is pure point. This will
be shown to happen if and only if the other is pure point as well
(read on for details and a discussion of related results.)

\smallskip
The motivation for our study comes from physics and, more precisely,
from the study of solids with long range aperiodic order and
crystal-like diffraction spectrum. Such solids are known as genuine
quasicrystals. The existence of quasicrystals is now a
well-established and widely accepted experimental fact. Even if
discussions about the precise structures will still be going on for a
while, the common feature of aperiodicity has opened a new chapter of
crystallography and solid state research.

The original discovery of quasicrystals \cite{SBGC,Ni} was somewhat
accidental and only possible through one of their most striking
features, namely their sharp Bragg diffraction with point symmetries
that are not possible for 3-dimensional crystals (such as $n$-fold
rotation axes with $n=8,10,12$, or icosahedral symmetry); see
\cite{MB} for a summary and \cite{BG} for a guide to the literature.
These experimental results called for a mathematical explanation and
created a subject now often referred to as mathematical diffraction
theory.

Mathematical diffraction theory deals with the Fourier
transform of the autocorrelation measure (or Patterson measure) of a
given translation bounded (possibly complex) measure $\omega$. Here, $
\omega$ is the mathematical idealisation of the physical structure of
a solid or, more generally, of any state of matter. In its simplest 
form, $\omega$ is just a Dirac comb, i.e., a countable 
collection of (possibly weighted) point measures which mimic
the positions of the atoms (and their scattering strengths). The
autocorrelation measure $\gamma^{}_\omega$ of $\omega$ (see below for
a precise definition) is then a positive definite measure. Its Fourier
transform $\widehat{\gamma^{}_\omega}$ is a positive measure, called
the diffraction measure, which models the outcome of a diffraction
experiment; see \cite{Cowley} for background material and details on
the physical justification of this approach.

Now, given this setting, one of the most obvious questions to address
is that for (the characterisation of all) examples of measures
$\omega$ with a diffraction measure $\widehat{\gamma^{}_\omega}$
that is pure point, i.e., consists of point measures only.

This was addressed by Bombieri and Taylor in \cite{BT}. However, no
rigorous answer could be given at that time. Soon after, Hof
\cite{Hof} showed that structures obtained from the cut-and-project
formalism \cite{KN} possess a pure point diffraction spectrum under
rather decent assumptions, and Solomyak started a rather systematic 
study of substitution dynamical systems with pure point spectrum 
in \cite{Boris2}. By now, large classes of examples are known
\cite{BM,LMS-1}, also beyond the class of ordinary projection
sets \cite{BMP}. Moreover, Schlottmann was able to free the
cut-and-project formalism from basically all specific properties of
Euclidean space \cite{Martin1} and established that all regular 
model sets are pure point diffractive \cite{Martin2}, see
also \cite{Moody} for a summary on model sets.

A cornerstone in many of these considerations is the use of ergodic
theory and the so-called Dworkin argument \cite{Dworkin} (see
\cite{Hof2,Martin2} as well). This argument links the diffraction 
spectrum to the dynamical spectrum. It can be used to infer pure 
point diffraction spectrum from pure point dynamical spectrum.  
These investigations heavily depend on the underlying point sets being 
point sets of finite local complexity (FLC). This, however, is not
necessary, as becomes clear from two alternative approaches, one for
general Dirac combs and measures on the basis of almost periodicity by
Moody and one of the authors \cite{BM}, and the other for so-called
deformed model sets by Bernuau and Duneau \cite{Duneau}.

Thus, at the moment, there is a considerable gap between the cases
that can be treated by the method of almost periodicity of measures
\cite{GdeL,BM} and those using ergodic theory and requiring FLC
together with unique ergodicity. It is the primary aim of this article
to narrow this gap.  This will be achieved by thoroughly analysing the
link between the diffraction spectrum and the dynamical spectrum given
by the Dworkin argument.

The analysis carried out below will also be a crucial ingredient in a
forthcoming paper of ours \cite{BL} which investigates the stability of
pure point diffraction. Namely, we will set up a perturbation theory
for pure point diffraction by studying deformations of dynamical
systems with pure point dynamical spectrum. Particular emphasis will
be put on deformed model sets and isospectral deformations of Delone
sets. Note that the deformation of model sets almost immediately leads
to point sets which violate FLC.

\smallskip
Let us now discuss our results in more detail.  The first step in our
approach is to choose a setting of measures rather than point
sets. Defining appropriate dynamical systems with measures on a
locally compact Abelian group will free us from essentially all
restrictions mentioned.  This setting is presented  in
Section~\ref{Measure}, where also the relevant topological questions
are discussed. The relationship between our measure dynamical systems
and point dynamical systems is investigated in Section~\ref{Delone}.
It is shown that the measure dynamical systems enclose the usual
Delone dynamical systems. This section introduces a topology on Delone
sets (and actually all closed subsets of the group) with very nice
compactness properties. The results generalise and strengthen the
corresponding considerations of \cite{Martin2,LS} and may be of
independent interest in further studies of point sets not satisfying
FLC.  The extension of diffraction theory and the Dworkin argument
(as developed for point sets in, e.g., \cite{Dworkin,Martin2}) to our
setting is achieved in Section~\ref{Diffraction}.

\smallskip
The main result of our paper is summarised in Theorem~\ref{Characterization} 
in Section~\ref{Main}. It states that, under some rather mild
assumptions,
\begin{itemize}
\item pure point dynamical spectrum is equivalent to pure point 
      diffraction spectrum.
\end{itemize} 
This generalises the main results of
Lee, Moody and Solomyak \cite{LMS-1} in at least three ways: It is not
restricted to dynamical systems arising from Delone sets, and, in fact,
not even to dynamical systems arising from point sets. It does not
need any condition of finite local complexity. It does not need an
ergodicity assumption for the invariant measure involved.

Let us mention that a generalisation of the main result of
\cite{LMS-1} sharing the last two features had already been announced
by Gou\'{e}r\'{e} \cite{Gouere-1,Gouere-3}, within the framework of
point processes and Palm measures, see his recent work \cite{Gouere-1}
for a study of this framework as well.  His result applies to the
point dynamical systems studied below in Section~\ref{Delone}.  Thus,
there is some overlap between his result an ours. However, in general,
our setting, methods and results are quite different from his, as we
leave the scenario of point sets. In fact, the measure theoretical
setting seems very adequate and natural in view of the physical
applications, where point sets are only a somewhat crude approximation
of the arrangements of scatterers.

Let us also mention that, in general, the diffraction spectrum and the 
dynamical spectrum can be of different type, as has been investigated 
by van Enter and Mi\c{e}kisz in \cite{EM}.

Our proof of the equivalence of the two notions of pure pointedness
relies on two results which are of interest in their own right. These
results are
\begin{itemize}
\item an abstract characterisation of pure point dynamical spectrum
  for arbitrary topological dynamical systems,
\item a precise interpretation of the Dworkin argument. 
\end{itemize}
Here, the abstract characterisation is achieved in 
Theorem~\ref{Abstractcriterion} in Section~\ref{Criterion}. Roughly 
speaking,
it states that a system has pure point spectrum once it has a lot of
point spectrum. The precise interpretation of the Dworkin argument is
given in Theorem~\ref{Spectralmeasure} in Section~\ref{Dynamical}. It
says that the diffraction measure is a spectral measure for a suitable
subrepresentation of the translation action at hand.

\smallskip
The relationship between suitable  subrepresentation with pure point
spectrum and the original representation can actually be analysed in
more detail.  To do so, we take a second look at the abstract theory in
Section  \ref{Further}. Namely, we discuss how the group of all eigenvalues
and the continuity of eigenfunctions is already determined by the set of 
eigenvalues and continuity of eigenfunctions associated to a suitable
subrepresentation with pure point spectrum. 

This material is rather general and may be of independent interest.
Here, we apply it  to our topological measure dynamical systems. This
gives a criterion for the continuity of the eigenfunctions. More importantly, 
it shows that the group of all eigenvalues is generated by the support
of the diffraction spectrum.  The validity of such a result was brought to 
our attention by R.~V.~Moody for  the case of point dynamical systems 
satisfying FLC \cite{Moody2002}. 

The material presented above, and the abstract strategy to prove our main 
result, can be adopted to study a measurable framework (as opposed to a
topological one). This will be analysed in the future.

\section{An abstract criterion}\label{Criterion}

In this section, we introduce some notation and provide a simple
result which lies at the heart of our considerations.  It is rather
general and might also be useful elsewhere.

\smallskip
Let $\varOmega$ be a compact topological space (by which we mean to
include the Hausdorff property) and $G$ be a locally
compact Abelian (LCA) group which is $\sigma$-compact. Let
\begin{equation} \label{def:action}
  \alpha \!:\, G \times \varOmega \;\longrightarrow\; \varOmega
\end{equation} 
be a continuous action of $G$ on $\varOmega$, where, of course,
$G\times \varOmega$ carries the product topology (later on, we will
specify it via $\alpha^{}_{t} (P) = t +P$ for $P\subset G$ and $t\in G$). 
Then, $\oO$ is called a {\em topological dynamical system}. 
The set of continuous functions on
$\varOmega$ will be denoted by $C(\varOmega)$. Let $m$ be a
$G$-invariant probability measure on $\varOmega$ and denote the 
corresponding set of square integrable functions on $\varOmega$
by $L^2(\varOmega,m)$. This space is equipped with the inner product 
$\langle f, g\rangle := \int \overline{f (\omega)} g(\omega)\dd m(\omega)$. 
The action $\alpha$ induces a unitary representation $T$ of $G$ on $L^2
(\varOmega,m)$ in the obvious way, where $T^t h$ is defined by 
$(T^t h)(\omega) := h(\alpha^{}_{-t} (\omega))$.  
Whenever we want to emphasise the dependence of the inner product and 
the unitary representation on the chosen invariant measure $m$, we write 
$\langle f, g\rangle_m$ and $T_m$ instead of $\langle f, g\rangle$ and $T$.

The dual group of $G$ is
denoted by $\widehat{G}$, and the pairing between a character
$\hat{s}\in \widehat{G}$ and an element $t\in G$ is written as
$(\hat{s}, t)$, which, of course, is a number on the unit circle,
compare \cite[Ch.~4]{Reiter} for background material.

A non-zero $h\in L^2 (\varOmega,m)$ is called an
\textit{eigenvector}\/ (or eigenfunction) of $T$ if there exists an
$\hat{s}\in \widehat{G}$ with $T^t h = (\hat{s}, t) h$ for every $t\in
G$. The closure  (in $L^2 (\varOmega,m)$) of the linear span of all 
eigenfunctions of $T$ will
be denoted by $\cHp$.

\smallskip
The following is a variant (and  an extension) of a
result from \cite{LMS-1}.
\begin{lemma}   \label{Algebra}
    Let\/ $\oO$ be a topological dynamical system with an invariant
    measure\/ $m$.  Then, $\cHp\cap C(\varOmega)$ is a subalgebra of\/
    $C(\varOmega)$ which is closed under complex conjugation and
    contains all constant functions. Similarly, $\cHp \cap
    L^\infty(\varOmega,m)$ is a subalgebra of\/ $L^\infty (\varOmega,m)$
    that is closed under complex conjugation and contains all constant
    functions. 
\end{lemma}
\begin{proof}  
We only show the statement about $\cHp\cap C(\varOmega)$. 
The other result can be shown in the same way.

The set $\cHp\cap C(\varOmega)$ is a vector
space because it is the intersection of two vector spaces. Moreover,
every constant (non-vanishing) function is obviously continuous and an
eigenvector of $T$ (with eigenvalue $1$, i.e., with the trivial
character $(\hat{s},t)\equiv 1$).

\smallskip
It remains to be shown that $\cHp\cap C(\varOmega)$ is closed under
complex conjugation and under forming products.

\smallskip
Closedness under complex conjugation: Let $f$ be an eigenfunction of
$T$ to, say, $\hat{s}$. Then, $\overline{f}$ is an eigenfunction of
$T$ to the character $\hat{s}^{-1}$. Here, of course, the inverse
$\hat{s}^{-1}$ of $\hat{s} \in \widehat{G}$ is given by $t\mapsto
\overline{ (\hat{s}, t)}$, where the bar denotes complex conjugation.
Using this, it is not hard to see that $\cHp$ is closed under complex
conjugation. As this is true of $C(\varOmega)$ as well, we see that
the intersection $\cHp\cap C(\varOmega)$ is closed under complex
conjugation.

\smallskip
Closedness under products: This is shown in Lemma 3.7 in \cite{LMS-1}
in the case that $m$ is not only translation invariant but also ergodic. 
To adopt their argument to the case at hand, we note that every
eigenfunction can be approximated arbitrarily well (in
$L^2(\varOmega,m)$) by bounded eigenfunctions via a simple cut-off
procedure. More precisely, if $f$ is an eigenfunction, then $|f|$ is
an $\alpha$-invariant function. Therefore, for an arbitrary $N>0$, 
the function
\begin{equation} \label{cutoff}
   f^N (\omega) \; := \; \begin{cases} f(\omega), 
           & |f(\omega)|\le N \\ 0 \, , 
           & \text{otherwise}
           \end{cases}
\end{equation}
is again an eigenfunction (with the same $\hat{s}$ as
$f$). Apparently, the $f^N$ converge to $f$ in $L^2(\varOmega,m)$ as
$N\rightarrow\infty$.

After this preliminary consideration, we can conclude the proof
following \cite{LMS-1}: Let two functions $f,g\in \cHp\cap C(\varOmega)$ 
be given. Then, $f g$ belongs to $C(\varOmega)$. It remains to be shown
that it belongs to $\cHp$ as well.

Choose $\varepsilon>0$ arbitrarily. Observe that $\|g\|_{\infty} < \infty$,
as $g\in C(\varOmega)$ with $\varOmega$ compact. 
Since $f$ is in $\cHp$, there exists
a finite linear combination $f'= \sum a_i f_i$ of eigenfunctions of
$T$ with
\begin{equation*}
   \|f - \sum a_i f_i\|_2 \; \le \; \frac{ \varepsilon}{\|g\|_\infty}\, .
\end{equation*} 
By the preliminary consideration around Eq.~(\ref{cutoff}), we can 
assume that all $f_i$ are
bounded functions. Thus, in particular, $\|f'\|_\infty < \infty$.

Similarly, choose another finite linear combination $g' = \sum b_j
g_j$ of bounded functions $g_j$ in $\cHp$ with
\begin{equation*}
   \|g - \sum b_j g_j\|_2  \; \le \;  
   \frac{\varepsilon}{\|f'\|_\infty}  \, .
\end{equation*}
Then, 
\begin{equation*}
   \| f g - f' g'\|_2 \; \le \; \|f'\|_\infty \, \| g - g'\|_2 + 
   \|g \|_\infty \, \|f - f'\|_2  \; \le \; 2\hspace{0.5pt} \varepsilon.
\end{equation*}
The proof is complete by observing that $f' g'$ is in $\cHp$ because
the product of bounded eigenfunctions is again a bounded
eigenfunction.  
\end{proof}

With Lemma~\ref{Algebra}, the following result is a rather direct 
consequence of the Stone-Weierstra{\ss} Theorem.
\begin{theorem}\label{Abstractcriterion} 
  Let\/ $\oO$ be a topological dynamical system with invariant probability
  measure\/ $m$. Then, the following assertions are equivalent.
\begin{itemize}
\item[\rm (a)] $T$ has pure point spectrum, i.e.,  
               $\cHp = L^2 (\varOmega,m)$.
\item[\rm (b)] There exists a subspace\/ $\cV\subset \cHp\cap C(\varOmega)$ 
               which separates points. 
\end{itemize}
\end{theorem}
\begin{proof}
(a) $\Longrightarrow$ (b): This is clear, as
one can take $\cV=\cHp\cap C(\varOmega)=C(\varOmega)$. \\ 
(b) $\Longrightarrow$
(a): By Lemma~\ref{Algebra}, $\cHp \cap C(\varOmega)$ is an algebra
which is closed under complex conjugation and contains the constant
functions. It also separates points as it contains a subspace, $\cV$,
with this property by (b). Thus, we can apply the Stone-Weierstra{\ss} 
Theorem (compare \cite[Thm.~4.3.4]{Ped}), to conclude that 
$\cHp \cap C(\varOmega)$ is dense in $C(\varOmega)$. By standard 
measure theory, see \cite[Prop.~6.4.11]{Ped}, $\cHp$ is then dense 
in $L^2(\varOmega,m)$ as well. As $\cHp$ is closed, statement (a)
follows.  
\end{proof}

\section{Measure dynamical systems} \label{Measure}

For the remainder of the paper, let $G$ be a fixed $\sigma$-compact
LCA group with identity $0$. Integration with respect to Haar measure
is denoted by $\int_G \ldots\dd t$, and the measure of a subset $D$ of
$G$ is denoted by $|D|$. The vector space of complex valued continuous 
functions on $G$ with compact support is denoted by $C_c (G)$. It is 
made into a locally convex space by the inductive limit topology, as 
induced by the canonical embeddings
\begin{equation*}
  C^{}_K (G) \; \hookrightarrow \; C_c (G)
  \; , \quad K\subset G \mbox{ compact}.
\end{equation*}
Here, $C_K (G)$ is the space of complex valued continuous functions on 
$G$ with support in $K$, which is equipped with the usual supremum norm
$\|.\|_{\infty}$. The support of $\varphi\in C_c(G)$ is denoted by 
$\supp (\varphi)$.  

The dual $C_c(G)^\ast$ of the locally convex space $C_c (G)$ is denoted
by $\MM$. The space $\MM$ carries the vague topology. 
This topology equals the weak-$\ast$ topology of
$C_c(G)^\ast$, i.e., it is the weakest topology which makes all
functionals $\mu\mapsto \mu(\varphi)$, $\varphi\in C_c (G)$,
continuous.  As is well known (see e.g. \cite[Thm.~6.5.6]{Ped} together 
with its proof), every $\mu \in \MM$ gives rise to a unique $\lvert\mu\rvert 
\in \MM$, called the {\em total variation}\/ of $\mu$, which satisfies
\begin{equation*}
 \lvert\mu\rvert(\varphi) \; = \; \sup\, \{ \lvert\mu (\psi) \rvert :
   \psi \in C_c (G,\RR) \mbox{ with } \lvert\psi\rvert\le\varphi \}
\end{equation*}
for every $\varphi \in C_c (G)_{+}$. 
Apparently, the total variation $\lvert\mu\rvert$ is positive, i.e.,
$\lvert\mu\rvert (\varphi) \ge 0$ for all $\varphi\in C_c (G)_{+}$.
In particular,  it   can be
identified with a measure on the $\sigma$-algebra of Borel sets of
$G$ that satisfies
\begin{itemize}
\item $\lvert\mu\rvert (K) < \infty$ for every compact set $K\subset G$,
\item $\lvert\mu\rvert (A) = \sup\, \{ \lvert\mu\rvert (K) : 
    K\subset A, K \mbox{ compact}\}$ for every Borel set $A\subset G$
\end{itemize}
(see, e.g., \cite[Thm.~6.3.4]{Ped}).
As $G$ is $\sigma$-compact, we furthermore have
\begin{itemize}
\item $\lvert\mu\rvert (A) = \inf \{ \lvert\mu\rvert (B) :
    A \subset B, B \mbox{ open} \}$ for every Borel set $A\subset G$
\end{itemize}
by \cite[Prop.~6.3.6]{Ped}, i.e., $\lvert\mu\rvert$ is an (unbounded)
regular Borel measure, in line with the Riesz-Markov Theorem
\cite[Thm.~IV.18]{RS}.
Finally, we note that there exists, by \cite[Thm.~6.5.6]{Ped},
a measurable function
$u\!: G\longrightarrow \CC$ with $\lvert u(t)\rvert = 1$ for
$\lvert\mu\rvert$-almost every $t\in G$ such that
\begin{equation} \label{polar}
  \mu(\varphi) \; = \; \int_G  \varphi\, u \dd \lvert\mu\rvert
  \quad \mbox{ for all } \varphi \in C_c (G) .
\end{equation}
This polar decomposition permits us to identify the elements of $\MM$ 
with the regular complex Borel measures on $G$, which is the Riesz-Markov
Theorem for this situation.

For later use, we also introduce some notation concerning Fourier
transforms and convolutions, compare \cite{Rudin,BF} for details.  
The Fourier transform of a quantity $q$
will always be denoted by $\widehat{q}$.  For $\varphi,\psi\in
C_c(G)$, we define the convolution $ \varphi\ast \psi$ by 
$\big(\varphi\ast\psi\big)(t) := \int_G \varphi(s) \psi(t-s)\dd s$ 
and the function
$\widetilde{\varphi}\in C_c(G)$ by $\widetilde{\varphi}(t):=
\overline{\varphi(-t)}$.
For $\mu\in \MM$ and $\varphi\in C_c (G)$, the convolution
 $\varphi\ast \mu$ is the function given by $\big(\varphi\ast\mu\big) (t)
 :=\int_G \varphi (t-s)\dd \mu(s)$.  
For two convolvable measures $\mu,\nu\in \MM$, the
 convolution $\mu \ast \nu$ is the element of $\MM$
 given by $\big(\mu\ast\nu\big) (\varphi) := \int_G \int_G \varphi(s + t)\dd
 \mu(s)\dd \nu(t)$ for $\varphi\in C_c(G)$; the measures
 $\widetilde{\mu}$ and $\overline{\mu}$ are defined by
$\widetilde{\mu}({\varphi}) :=
 \overline{\mu(\widetilde{\varphi})}$ and 
$\overline{\mu} (\varphi):=
 \overline{\mu(\overline{\varphi})}$, respectively.
For $\mu\in \MM$ and a
 measurable set $B\subset G$, we denote the restriction of $\mu$ to
 $B$ by $\mu^{}_B$. Finally, for $x\in G$, we define the measure
 $\delta_x$ to be the normalised point measure at $x$.

\smallskip
We will consider actions of $G$ on spaces consisting of measures on
$G$.  The relevant set of measures will be defined next.
\begin{definition}  \label{mess-def}
  Let\/ $C>0$ and a relatively compact open set\/ $V\!$ in\/ $G$ be given.
  A measure\/ $\mu\in \MM $ is called\/ $(C,V)$-{\em translation bounded} if\/ 
  $|\mu|(t + V) \le C$ for all\/ $t \in G$. It is simply called\/ {\em translation 
  bounded} if there exist\/ $C,V\!$ such that it is\/ $(C,V)$-translation 
  bounded. The set of all\/ $(C,V)$-translation bounded measures is denoted 
  by\/ $\MCV$ and the set of all translation bounded measures by\/ $\MTB$.
\end{definition}

The vague topology on $\MM$ has very nice features when restricted to the
translation bounded measures.
\begin{theorem}\label{Compact} 
  Let\/ $C>0$ and a relatively compact open set\/ 
  $V\!$ in\/ $G$ be given. Then, $\MCV$ is a compact Hausdorff space. 
  If\/ $G$ is second countable, $\MCV$ is metrisable.
\end{theorem}

To prove the theorem, we start with the following simple result from measure 
theory.
\begin{prop} \label{nuetzlich} 
  Let\/ $\mu\in\MM$ and a relatively compact open set\/ $V$ in\/
  $G$ be given. Then, $|\mu| (V) = \sup\,\{
  \lvert\mu(\varphi)\rvert : \varphi\in C_c (G) \mbox{ with }\,
  \supp(\varphi)\subset V \!\mbox{ and }\, \|\varphi\|_{\infty} \le 1 \}$.
\end{prop}
\begin{proof}  
Denote the supremum in the statement by $S$. Recall the polar decomposition and 
the definition of $u$ in  Eq.~(\ref{polar}). As $u$ is $\lvert\mu\rvert$ almost 
surely equal to  $1$, a direct calculation easily gives $S\le |\mu| (V)$.
Conversely, as $C_c (V)$ is dense in
$L^1 (V,\lvert\mu\rvert^{}_{V})$  (see, e.g., \cite[Prop.~6.4.11]{Ped}), there 
exists a sequence $(\varphi^{}_{n})$ in $C_c(V)$ with $\varphi^{}_{n} \, 
\xrightarrow{n\to\infty} \; \overline{u}\cdot 1^{}_{V}
  \; \mbox{in } L^1 (V,\lvert\mu\rvert^{}_{V})\,.$
By the $\lvert\mu\rvert$-almost sure boundedness of $u$, this implies
$ u \varphi^{}_{n} \xrightarrow{n\to\infty} u \overline{u}\cdot 1^{}_{V} =1^{}_{V} $
in $L^1 (V,\lvert\mu\rvert^{}_{V})$. Then, a short calculation,
invoking \eqref{polar} again, shows $ S \ge \lim_{n\to\infty} \;
\lvert\mu(\varphi^{}_{n})\rvert 
\;\, = \;\, \lvert\mu\rvert (V).$
This  proves the proposition.
\end{proof}

\noindent {\it Proof of Theorem~$\ref{Compact}$.}  
By definition of
$\MCV$, for each $\varphi \in C_c (G)$, there exists a radius
$R(\varphi)>0$ such that $\mu(\varphi)\in \overline{B^{}_{R(\varphi)}}$
for every $\mu\in \MCV$, where $B_r$ is the (open) ball of radius $r$
around $0$ and $\overline{B_r}$ its closure.
Thus, we can consider $\MCV$ as a subspace of
$\varPi:= \prod_{\varphi\in C_c (G)} \overline{B^{}_{R(\varphi)}}$
equipped with the product topology via the embedding
\begin{equation*} 
    j\! :\, \MCV \; \hookrightarrow \; \varPi
    \, , \quad \big(j(\mu)\big)(\varphi):= \mu(\varphi).
\end{equation*}
As $\varPi$ is obviously a compact Hausdorff space, this shows
immediately that $\MCV$ is relatively compact and Hausdorff. It
remains to be shown that $j(\MCV)$ is closed. This is a direct
consequence of Definition~\ref{mess-def} together with 
Proposition~\ref{nuetzlich}. 

The statement about metrisability is standard: if $G$ is second countable,
there exists a countable dense subset $\{\varphi_n : n\in \NN\}$ in
$C_c (G)$.  Then,
\begin{equation*} 
   d\! :\, \MCV \times \MCV\longrightarrow \RR \, , 
           \quad d(\mu,\nu) \; := \;
    \sum_{n\in \NN} 2^{-n} \frac{ |\mu(\varphi_n) - \nu(\varphi_n)|}
                   { 1 + |\mu(\varphi_n) - \nu(\varphi_n)|}
\end{equation*}
gives a metric on $\MCV$ which generates the topology.  
\hfill \qedsymbol

\smallskip
Having discussed the topology of $\MM$, we can now introduce the
topological dynamical systems associated to subsets of $\MM$. To do
so, we will use the obvious action $\alpha$ of $G$ on $\MM$ given by
\begin{equation*}
   \alpha \!:\, G\times \MM \longrightarrow \MM \, ,\quad 
        (t,\mu)\mapsto
   \alpha^{}_t(\mu)  := \delta_t \ast \mu, 
\end{equation*}
or, more explicitly,  $\big(\alpha^{}_t (\mu)\big) (\varphi) =  
\int_G \varphi(t + s)\dd \mu(s).$ We use the same symbol for the 
action as in Eq.~(\ref{def:action}), since misunderstandings are
unlikely, and we will usually write $\alpha^{}_t \mu$
for $\alpha^{}_t(\mu)$.

This action is compatible with the topological structure of $G$ and $\MM$.
\begin{prop} \label{continuity} 
  Let\/ $C>0$ and a relatively compact open set\/ $V\subset G$ be given.
  Then, the action\/ $\alpha \!:\, G\times
  \MCV \longrightarrow \MCV$ is continuous, where\/ $G\times \MCV$ 
  carries the product topology.
\end{prop}
\begin{proof}  
Let $(t_\iota,\mu_\iota)$ be a net in $G\times\MCV$ 
converging to $(t,\mu)$. We have to show that the net
$(\alpha^{}_{t_\iota}(\mu_\iota))$ converges to $\alpha^{}_t(\mu)$,
i.e., we have to check that
\begin{equation*} 
   \int \varphi(s + t_\iota)\dd \mu_\iota (s) \; \longrightarrow \;
   \int \varphi(s + t)\dd \mu(s) \, , \quad
   \mbox{for all $\varphi\in C_c (G)$}.
\end{equation*}
By $t_\iota\longrightarrow t$, there exists an index $\iota^{}_0$ and a
compact set $K$ such that the support of $\varphi$ and the supports of all 
$\varphi_\iota$ with $\iota\geq \iota^{}_0$ are contained in $K$. Moreover, 
as $\mu_\iota\in \MCV$ for every $\iota$, there exists a constant $C'$ with 
$|\mu|(K)\leq C'$ as well as $|\mu_\iota| (K)\leq C'$ for all $\iota$. Now, 
the desired statement follows easily from 
\begin{eqnarray*}
  |\mu_\iota (\varphi (. + t_\iota)) - \mu  
           (\varphi(. + t))|
  &\leq & |\mu_\iota (\varphi (. + t_\iota)) - \mu_\iota 
  (\varphi( . + t))|  +|\mu_\iota (\varphi (. + t)) - \mu  
  (\varphi( . + t))|\\
  &\leq & |\mu_\iota| (C) \| \varphi (. + t_\iota) - \varphi
     ( . + t)\|_\infty +|\mu_\iota (\varphi (. + t)) - \mu  
     (\varphi( . + t))|,
\end{eqnarray*}
as both terms on the right hand side tend to zero for  
$t_\iota \longrightarrow t$ and $\mu_\iota \longrightarrow \mu$. 
\end{proof}

The dynamical systems we are interested in are defined as follows. 
\begin{definition} 
  $\oO$ is called a\/ {\em dynamical system} on the translation bounded
  measures on\/ $G$\/ {\rm (TMDS)} if there exist a constant\/ $C >0$ and a
  relatively compact open\/ $V\subset G$ such that\/ $\varOmega$ is a
  closed subset of\/ $\MCV$ that is invariant under the\/ $G$-action\/
  $\alpha$.
\end{definition}

\smallskip
\noindent{\bf Remarks.}  
(a) The $\alpha$ invariant subsets of $\MM$
are called \textit{translation stable sets} in \cite{GdeL}. Thus, a
TDMS is just a closed translation stable subset of $\MCV$. \newline
(b) The space $\varOmega$ of a TMDS is always compact by
Theorem~\ref{Compact} and the action $\alpha$ is continuous by
Proposition~\ref{continuity}. Thus, a TMDS is a topological
dynamical system in the sense of Section~\ref{Criterion}. \newline (c)
The considerations of this section (and those of the next) do not use
commutativity of the underlying group $G$. They immediately extend to
arbitrary $\sigma$-compact locally compact groups. But since we need
harmonic analysis later on, we stick to Abelian groups here.

\section{Point dynamical systems}\label{Delone}

This section has two aims. Firstly, we present an abstract topological
framework which allows us to treat point dynamical systems which are
not of finite local complexity. Secondly, we
show how these systems fit into our setting of measure dynamical
systems. As for the first aim, we actually introduce a topology on the
set of all closed subsets of $G$. For the case of $\RR^d$, this
topology has already been studied by Stollmann and one of the authors
in \cite{LS}. Our extension to arbitrary locally compact groups is
strongly influenced by the investigation of Schlottmann \cite{Martin2}
(which, however, is restricted to FLC systems). 

\smallskip

As pointed out by the referee, the topology we introduce can
also be obtained as a special case of a topology
introduced by Fell in \cite{Fell}. This is further discussed in the
appendix.  Given this connection, Theorem \ref{Topology} below is
a corollary of Theorem 1 in \cite{Fell}. For this reason, we only 
give an outline of how it can be established in our setting.

\smallskip

We start by defining the relevant sets of points. 
\begin{definition} 
Let\/ $G$ be a\/ $\sigma$-compact\/ {\rm LCA} group, and\/ $V\!$ an 
open neighbourhood of\/ $0$ in\/ $G$.
\begin{itemize}
\item[\rm (a)]  A subset\/ $\varLambda$ of\/ $G$ is called\/ $V\!$-{\em discrete} 
   if every translate of\/ $V\!$ contains at most one point of\/ $\varLambda$. 
   The set of all\/ $V\!$-discrete subsets of\/ $G$ is denoted by\/ $\cDV$. 
\item[\rm (b)] A subset of\/ $G$ is called\/ {\em uniformly discrete} if it is\/ 
   $V\!$-discrete for some\/ $V\!$. The set 
   of all uniformly discrete subsets of\/ $G$ is denoted as\/ $\cUD$.  
\item[\rm (c)] The set of all discrete and closed subsets of\/ $G$ will
   be denoted by\/ $\mathcal{D} (G)$. 
\item[\rm (d)] The set of all closed subsets of\/ $G$ is denoted as\/ $\cC$.
\item[\rm (e)] A subset\/ $\varLambda$ of\/ $G$ is called\/ {\em relatively dense} 
  if there is a compact\/ $K$ with\/ $ \varLambda + K = G$. 
\end{itemize}
\end{definition}

Note that a uniformly discrete subset of $G$ is closed. As it
presents no extra difficulty, we will actually topologise not only
$\mathcal{D}(G)$ but rather the larger set $\cC$. This will be done by
providing a suitable uniformity (see \cite[Ch.~6]{Kel} for details on
uniformities). Namely, for $K\subset G$ compact and $V\!$ a
neighbourhood of $0$ in $G$, we set
\begin{equation*}
   \UKV \; := \; \{(P_1,P_2) \in \cC \times \cC : 
   P_1 \cap K \subset P_2 + V \mbox{ and } P_2 \cap K \subset P_1 + V\}. 
\end{equation*}
It is not hard to check that
\begin{equation*} 
   (P,P) \in \UKV,\;\: \UKV= \UKV^{-1}, \;\: U_{K_1\cup K_2, V_1\cap
   V_2} \subset U_{K_1,V_1} \cap U_{K_2,V_2},\:\; U_{K-W, W} \circ
   U_{K-W,W} \subset U_{K,V}
\end{equation*} 
for $V\!$ a neighbourhood of $0$, $W$ a compact neighbourhood of
$0$ with $W + W \subset V$, $K$ in $G$ compact, and $P$ any closed
subset of $G$. Here, on sets $U,U_1,U_2$ consisting of ordered pairs,
we define $U^{-1} := \{ (y,x) : (x,y)\in U\}$ and
\begin{equation*} 
   U_1 \circ U_2 \; := \; \{(x,z) :\, \exists\, y\in G\,
   \mbox{ with }\, (x,y)\in U_1\;\: 
   \mbox{and}\;\: (y,z) \in U_2\}. 
\end{equation*}
This guarantees that $\{\UKV : K\,\mbox{compact}, \,V \,\mbox{open with 
$ 0\in V$}\}$ generates a uniformity, and hence a topology on $\cC$ via 
the neighbourhoods
\begin{equation*}
    \UKV (P)  \; := \; \{ Q : (Q,P)\in \UKV\}\, ,
    \quad P\in \cC\, .
\end{equation*}
Note that we could equally well generate the same uniformity with
$V$ running through compact neighbourhoods of $0\in G$.
\begin{definition}
   The topology defined this way is called the\/ {\em local rubber topology}
   {\rm (LRT)}.
\end{definition}

This topology essentially means that two sets $P_1, P_2$ are
close if they ``almost'' agree on large compact sets.

\smallskip
Fundamental properties of the LRT are given in the following result (see
\cite{LS} for an earlier result on $\RR^d$).
\begin{theorem} \label{Topology} 
  With the {\rm LRT}, the set\/ $\cC$ of closed subsets of\/ $G$ 
  is a compact Hausdorff space.  If the topology of\/  $G$ has a countable base, 
  then\/ $\cC$ is metrisable. 
\end{theorem}
\begin{proof} 
The set $\cC$ is Hausdorff, as the intersection of all $\UKV$ 
contains only the diagonal set $\{(P,P) : P \in \cC \}$, see
\cite[Ch.~6]{Kel}.

\smallskip
We next show completeness: Let $(P_\iota)_{\iota\in I}$ be a Cauchy
net in $\cC$, where $I$ is an index set directed by $\le$, compare
\cite[Ch.~2]{Kel}. 

We have to  show that the Cauchy net
converges to a closed subset of $G$, hence an element of $\cC$.
To this end, we introduce the set $P$ of those $x\in G$ such that, 
for every neighbourhood $V$ of $0$, there exists 
$\iota^{}_{x,V}\in I$ with
\begin{equation} \label{define-P}
   (x + V) \cap P_{\iota}\; \neq\; 
   \varnothing \quad \mbox{for all $\iota \ge \iota^{}_{x,V}$}.
\end{equation}
It is not hard to see that $P$ is closed.  $P$ will turn out to be  
the limit of our Cauchy net.

\smallskip
To show this,  let a compact $K$ in $G$ and a neighbourhood $V$ of $0$ be
given.  We have to provide an $\iota^{}_{K,V}$ with
\begin{equation*}
  P\cap K \, \subset \, P_{\iota} + V \quad\mbox{and}\quad 
  P_\iota \cap K \, \subset \,  P + V
\end{equation*}
for every $\iota \ge\iota^{}_{K,V}$.

\smallskip
Rather direct arguments show existence  of $\iota^{}_{K,V}$ with  
$P\cap K \subset P_{\iota} + V$  for all $\iota \ge \iota^{}_{K,V}$.
We next establish the other inclusion.  By a compactness argument, 
\begin{equation}
\label{claim}
   P\cap C \; \neq \; \varnothing,
\end{equation}
 whenever $C\subset G$ is compact
and $P_{\kappa}\cap C\neq\varnothing$ for all\/ $\kappa\ge\kappa^{}_0$
and some\/ $\kappa^{}_0\in I$.
Assume, without loss of generality, that $V$ is compact and symmetric.
As $(P_\iota)$ is a Cauchy net, there exists a $\iota^{}_{K,V}$ with
$(P_\iota,P_\kappa) \in U_{K,V}$ for all $\iota,\kappa \ge
\iota^{}_{K,V}$.  Consider $\iota \ge \iota^{}_{K,V}$ and choose an
arbitrary $q\in P_\iota \cap K$. Then,
\begin{equation*}
 (q + V) \cap P_\kappa \; \neq \; \varnothing
\end{equation*}
for every $\kappa \ge \iota^{}_{K,V}$, as $V$ is symmetric and
$P_\iota \cap K \subset P_\kappa + V$ by $(P_\iota,P_\kappa) \in
U_{K,V}$.  Thus, \eqref{claim} gives existence of a $p\in P$ with
$p\in q + V$. In particular, invoking the symmetry of $V$ once more,
we have $q \in p + V \subset P + V$. As $q\in P_\iota \cap K$ was
arbitrary, we infer the desired inclusion $ P_\iota \cap K \;
\subset\; P + V$.

These  considerations  prove the desired completeness statement.  

\smallskip
Finally, we show compactness of $\cC$. As $\cC$ is complete, it
suffices to prove it is precompact.  Thus, for any given $K$ compact
and $V\!$ an open neighbourhood of $0$ in $G$, we have to provide a
natural number $n$ and $P_i\in \cC$, $1 \le i \le n$, such that
$   \cC \; \subset \; \bigcup_{i=1}^n \UKV (P_i).$
Since $\UKV (P) \supset U_{K,V\cap (-V)} (P)$ for all $P\in \cC$ and
$V\cap (-V)$ is symmetric, we can assume, without loss of generality,
that $V\!$ is symmetric (i.e., $V= -V$). As $K$ is compact, there
exists a finite set $D\subset K$ with
$   K\subset D + V.$  Direct calculations then give 
\[
   \cC \; \subset \; \bigcup_{i=1}^n\UKV (D_i),
\] 
where $D_i$, $1\le i \le n$, is an enumeration of the power set of $D$. 

\smallskip
The statement about metrisability is a direct consequence of 
\cite[Thm.~6.13]{Kel} and the remark thereafter. 
\end{proof}

\smallskip
Having topologised $\cC$, and thus $\cUD$ as well, we can now
introduce our point dynamical systems. The natural action of $G$ on
$\cC$ by translation will also be denoted by $\alpha$. Explicitly, we
define $\alpha^{}_{t} (P) = t + P$ for $P\in\cC$, where $t + P = \{
t+x : x\in P\}$ as usual.
\begin{definition} 
Let\/ $\varOmega$ be a subset of\/ $G$ and\/ $\alpha$
the translation action just defined.
\begin{itemize}
\item[\rm (a)] The pair\/ $\oO$ is called a\/
   {\em  set dynamical system} if\/ 
    $\varOmega$ is a closed subset of\/ $\cC$ which is invariant under\/ 
    $\alpha$. 
\item[\rm (b)] A set dynamical system\/ $\oO$ is called a\/
    {\em point dynamical system} 
    if\/ $\varOmega$ is a subset of\/ $\cDV$ for some open 
    neighbourhood\/ $V\!$ of\/ $0$. 
\item[\rm (c)] A point dynamical system\/ $\oO$ is called a\/
    {\em Delone dynamical system}, 
    if every element of\/ $\varOmega$ is a relatively dense subset 
    of\/ $G$.
\end{itemize}
\end{definition}

It follows from Theorem~\ref{Topology} that a set dynamical system is
indeed a topological dynamical system in the sense of
Section~\ref{Criterion}.

 \smallskip
A special way of obtaining set dynamical systems is the following:
Choose $P\in \cC$. Then, the LRT-closure $X(P)$ of the orbit
$\{\alpha^{}_t(P) : t\in G\}$ of $P$ in $ \cC$ is a closed
$\alpha$-invariant subset of $\cC$, hence compact.  Thus,
$(X(P),\alpha)$ is a set dynamical system. We say that $P$ and $P'$
are {\em locally indistinguishable}\/ if $P\subset X(P')$ and $P' \subset
X(P)$, hence if $X(P) = X(P')$. The set of all $P'$ which are locally
indistinguishable form $P$ is called the RLI class of $P$, written as
RLI$(P)$. Here, as before, the letter R stands for ``rubber''.

\smallskip
Given these notions, we can characterise minimality of a set dynamical system 
in the following way, which extends \cite[Prop.~3.1]{Martin2} to our setting.
\begin{prop} 
   Let\/ $P\in \cC$ be given. Then, the following 
   assertions are equivalent.
\begin{itemize}
\item[\rm (a)] The set dynamical system\/ $(X(P), \alpha)$ is minimal.
\item[\rm (b)] $X(P) = {\rm RLI}(P)$.
\item[\rm (c)] The set\/ $\{ t\in G : (t+ P,P) \in \UKV\}$ is relatively 
    dense in\/ $G$ for all compact\/ $K\subset G$ and all open 
    neighbourhoods\/ $V\!$ of\/ $0$. 
\item[\rm (d)] The set\/ $P$ is repetitive, i.e., for every $K$ compact 
    and every open neighbourhood $V\!$ of\/ $0$, there is a compact set\/ 
    $C= C(K,V)\subset G$ such that, for\/ $t_1,t_2\in G$, there is an\/ 
    $s\in C$ with\/ $(t_1 + P, s + t_2 + P)\in \UKV$. 
\end{itemize}
\end{prop}
\begin{proof}  The equivalence of (a) and (b) is clear by the
definition of minimality and the local indistinguishability class. 
This is also known as Gottschalk's Theorem, compare \cite[Thm.~4.1.2]{Pet}.
The remaining assertions can be shown by mimicking and slightly
extending the proof of \cite[Prop.~3.1]{Martin2}.  
Since we do not use these results later on, we skip further details.
\end{proof}

\smallskip
Having discussed point dynamical systems, we can now relate them to
special dynamical systems on measures. The connection relies on the map
\begin{equation*} 
   \delta \!:\, \cUD \, \longrightarrow \, \MTB \, , \quad 
   \delta (\varLambda) := \sum_{x\in \varLambda} \delta_x  \, ,
\end{equation*} 
Note that $\delta$ is indeed a map into $\MTB$, as every $\varLambda\in
\cUD$ is uniformly discrete.
\begin{lemma}\label{Eigenschaftendelta}
 The map\/ $\delta \!:\, \cUD\longrightarrow \MTB$ is injective, continuous
 and compatible with the action of\/ $G$. The inverse\/ $\delta^{-1} \!:\,
 \delta(\cUD) \longrightarrow \cUD$ is continuous as well.
\end{lemma}
\begin{proof}  Injectivity and compatibility with the action
of $G$ are immediate. In particular, using $\delta_{t+x} =
\delta_{t}\ast\delta_{x}$, one can check that
\begin{equation*}
   \delta(\alpha^{}_{t} (\varLambda)) \, = \,
   \delta(t+\varLambda) \, = \,
   \delta^{}_{t} \ast \delta(\varLambda) \, = \,
   \alpha^{}_{t} (\delta(\varLambda))\, .
\end{equation*} 
To show continuity of $\delta$, we have to show
$\delta(\varLambda_\iota)(\varphi) \longrightarrow
\delta(\varLambda)(\varphi)$ for all $\varphi\in C_c (G)$, whenever
$\varLambda_\iota\longrightarrow \varLambda$. Let $V\!$ be an open
neighbourhood of $0$ in $G$ such that $\varLambda\in \cDV$. Let
$\varphi\in C_c (G)$ be given, so $\supp(\varphi)$ is compact. By a
simple partition of unity argument, we can now write
$ \varphi= \sum_{j=1}^n \varphi_j$, where each $\varphi_j$,
$j=1,\ldots, n$, has its support in a set of the form $W + t_j$, with $W
+ W \subset V$. Now, for such $ \varphi_j$, the convergence
$\delta(\varLambda_\iota)(\varphi_j) \longrightarrow
\delta(\varLambda)(\varphi_j)$ follows easily. This yields the desired
convergence for $\varphi$.  The continuity of $\delta^{-1}$ can be shown
similarly. 
\end{proof}

\smallskip
As can be seen from simple examples, $\delta(\cUD)$ is, in general, not
closed in $\MTB$ and the LRT is not the same as the vague topology on
$\delta(\cUD)$. For example, if $(x_n)$ is a sequence in $G$ with $x_n
\longrightarrow 0$ and $x_n \neq 0$, then, $\varLambda_n :=
\{0,x_n\}\in \cUD$ with $\varLambda_n \longrightarrow \{0\}$ in the LRT.
However, $\delta(\varLambda_n)\longrightarrow 2 \delta_0$ in the vague
topology, and $2\delta_0$ does not belong to $\delta (\cUD)$. Nevertheless, 
the following still holds.
\begin{prop} 
  Let\/ $V\!$ be an open neighbourhood of\/ $0$ in\/ $G$. Then, $\cDV$ is
  compact in the\/ {\rm LRT} and\/ $\delta(\cDV)$ is compact in the vague
  topology.
\end{prop}
\begin{proof}  
By compactness of $\cC$ and continuity of
$\delta$, it suffices to show that $\cDV$ is closed.  Since $V\!$ is
open, this is easy.  
\end{proof}

\smallskip
Our main result on the relationship of point dynamical systems and the 
framework of measure dynamical systems reads as follows. 
\begin{theorem} \label{Beziehung} 
  If\/ $\oO$ is a point dynamical system, the restriction\/
  $\delta|_\varOmega \!:\, \varOmega \longrightarrow \delta(\varOmega)$
  of\/ $\delta$ to\/ $\varOmega$ establishes a topological conjugacy between
  the dynamical systems\/ $\oO$ and\/ $(\delta(\varOmega),\alpha)$.
\end{theorem}
\begin{proof}
By Lemma~\ref{Eigenschaftendelta}, $\delta|_\varOmega$ is an
injective and continuous 
map which is compatible with the group action. As $\varOmega$
is compact, so is $\delta(\varOmega)$. By a standard argument
\cite[Prop.~1.6.8]{Ped}, $\delta|_\varOmega$ is then a homeomorphism
between $\oO$ and $(\delta(\varOmega),\alpha)$. Together, this
establishes the topological conjugacy.
\end{proof}

\smallskip
We finish this section by briefly discussing the relationship of the LRT
and the topology usually considered for Delone dynamical systems with
the FLC property.  A thorough discussion of the latter topology
has been given in \cite{Martin2}. This discussion actually gives a
topology on the closed subsets of $G$ (though this is not explicitly
noted in \cite{Martin2}). This topology will be called the local matching
topology (LMT).

The  definition  of the LMT in  \cite{Martin2} shows immediately that 
the LRT is coarser than the LMT. Thus, the identity
\begin{equation*}
   id \!:\, (\cC, {\rm LMT})\, \longrightarrow \, (\cC, {\rm LRT})
      \, , \quad P \mapsto P
\end{equation*}
is continuous. This yields the following result, which essentially
shows that our way of topologising the uniformly discrete sets
coincides with the usual topology when restricted to sets of finite
local complexity.
\begin{prop} 
    Let\/ $\varOmega$ be a subset of\/ $\cC$. If\/ $\varOmega$ is
    compact in the\/ {\rm LMT}, then\/ $\varOmega$ is compact in the\/
    {\rm LRT} as well, and the two topologies agree on\/ $\varOmega$.
\end{prop}
\begin{proof}  
The restriction $id_\varOmega \!:\, (\varOmega,
{\rm LMT})\longrightarrow (\varOmega, {\rm LRT})$ of the identity
to $\varOmega$ is continuous. Thus, as $(\varOmega, {\rm LMT})$ is
compact, so is its image $(\varOmega,{\rm LRT})$. Now, continuity of
the inverse is standard, cf.\ \cite[Prop.~1.6.8]{Ped}. 
Thus, the two topologies agree.
\end{proof}

\section{The diffraction spectrum}\label{Diffraction}
The basic concepts in the mathematical treatment of
diffraction experiments are the \textit{autocorrelation measure} and
the \textit{diffraction measure}. In the context of Delone dynamical
systems, these concepts have been developed and investigated  in a series
of articles by theoretical physicists and mathematicians
\cite{BH,BM,Dworkin,Hof,Hof2,Martin1,Martin2}. This will now be generalized 
and extended to our measure dynamical systems.  More
precisely, we show the existence of the autocorrelation measure by a
limiting procedure, provided certain ergodicity assumptions hold.

Let us mention that we will provide an alternative approach
to these quantities later on. It will be more general in that it does 
not need an ergodicity assumption.

\smallskip

To phrase our results, we need two more pieces of notation. Firstly,
recall from \cite{Martin2} that a sequence $(B_n)$ of compact subsets
of $G$ is called a \textit{van Hove sequence} if
\begin{equation*} 
   \lim_{n\to \infty} \frac{|\partial^K B_n|}{|B_n|} \; = \; 0
\end{equation*}
for all compact $K\subset G$. Here, for compact $B,K$, the 
``$K$-boundary''  $\partial^K B $  of $B$ is defined as
\begin{equation*} 
   \partial^K B \; := \; \overline{((B + K)\setminus B)} \cup 
    (\overline{G\setminus B} - K) \cap B,
\end{equation*}
where the bar denotes the closure. The existence of van Hove
sequences for all $\sigma$-compact LCA groups is shown in
\cite[p.~249]{Martin2}, see also Section 3.3 and Theorem (3.L)
of \cite[Appendix]{Tempelman}. Moreover, every van Hove sequence is
a F{\o}lner sequence, i.e., $\lvert B_n \triangle (B_n + K) \rvert /
\lvert B_n \rvert \xrightarrow{n\to\infty} 0$, for every compact set 
$K\subset G$, where $\triangle$ denotes the operator for the 
symmetric difference of two sets, compare Section 3.2 of
\cite[Appendix]{Tempelman}; for a partial converse, consult
Theorem (3.K) of the same reference.

Secondly, for $\varphi \in C_c (G)$ and $\mu \in \MM$, we define
\begin{equation*}
   f_\varphi (\mu) \; := \; \big(\varphi \ast \mu\big) (0) \; = \; 
   \int_G \varphi(-s)\dd \mu(s).
\end{equation*} 
This gives a way to ``push'' functions from $C_c (G)$ to
$C(\varOmega)$. In fact, up to the sign, $f_\varphi$ is just the
canonical embedding of $C_c (G)$ into its bidual $\MM^\ast$. Basic
features of the map $\varphi \mapsto f_{\varphi}$ are gathered in the
following lemma.
\begin{lemma}\label{con}
Recall the definition\/ $\alpha^{}_t \mu = \delta_t \ast \mu$ for
measures\/ $\mu$, and set\/ $\beta_t (\varphi) = \delta_t \ast \varphi$
for functions\/ $\varphi$. Then one has:
\begin{itemize}
\item[\rm (a)] The function\/ $f_\varphi \!:\, \MM \longrightarrow \CC$, 
   $\mu \mapsto f_\varphi (\mu)$ is continuous, for all\/ 
   $\varphi \in C_c (G)$.  
\item[\rm (b)] If\/ $\oO$ is a\/ {\rm TMDS}, the map\/ $f \!:\,
   C_c(G)\longrightarrow C(\varOmega),\:\;\varphi \mapsto f_\varphi$,
   is linear, continuous and compatible with the action of\/ $G$ in that\/
   $f_\varphi(\alpha^{}_t \mu)=f_{\beta_t (\varphi) }(\mu)$.
\end{itemize}
\end{lemma}
\begin{proof} 
(a) For $\varphi\in C_c (G)$, we have $f_\varphi (\mu)= \mu(\varphi_{\!\_})$, 
where $ \varphi_{\!\_}(t) = \varphi(-t)$. Thus, continuity of
$f_\varphi$ is immediate from the definition of the topology on $\MM$.

(b) Linearity of the map $f$ is obvious. To show continuity of $f$, recall 
that $C_c (G)$ is equipped with the inductive limit topology induced from 
the embeddings $C^{}_K (G) \hookrightarrow  C_c (G)$ with $K\subset G$ 
compact. Thus, it suffices to show the continuity of the map
\begin{equation*}
   f_K\! :\, C^{}_K (G) \; \hookrightarrow \; C_c (G) \, ,
   \quad \varphi \mapsto f_\varphi
\end{equation*}
for every  compact $K$ in $G$. 
So, let $(\varphi_\iota) $ be a net in $ C^{}_K (G)$  converging to 
$\varphi\in C^{}_K (G)$. Then,
$\|\varphi_\iota - \varphi\|_\infty \longrightarrow 0$, and 
$\supp(\varphi), \supp (\varphi_\iota) \subset K$ for all $\iota$. As 
$\varOmega \subset \MCV$
with suitable $C,V$, this easily implies $f_{\varphi_\iota}
\longrightarrow f_{\varphi}$.  Finally, a direct calculation shows
$f_\varphi(\alpha^{}_t \mu)=f_{\beta_t (\varphi)}(\mu)$.  
\end{proof}

\smallskip

\noindent{\bf Remark.}  The lemma is particularly interesting as there
does not seem to exist any canonical map from $\varOmega $ to $G$ or
from $G$ to $\varOmega$ in our setting (let alone a map which is
compatible with the corresponding group actions). However, if one
views a function $\varphi$ as the Radon-Nikodym derivative of a
measure that is absolutely continuous with respect to the Haar measure
of $G$, the action $\alpha^{}_{t}$ induces $\beta_t$ as defined.

\smallskip
Now, our result on existence of the autocorrelation function reads as
follows.
\begin{theorem} \label{Dworkin} 
Let\/ $\alpha$ be the translation action of\/ $G$ on\/ $\MTB$ as 
introduced above. 
\begin{itemize} 
\item[\rm (a)] If\/ $\oO$ is a uniquely ergodic\/ {\rm TMDS}, there
   exists a translation bounded measure\/ $\gamma$ on\/ $G$ such that the
   sequence\/ $(\frac{1}{|B_n |}\, \widetilde{ \omega_{B_n}} \ast
   \omega_{B_n})$ converges, in the vague topology, to\/ $\gamma$ for
   every van Hove sequence\/ $(B_n)$ and every\/ $\omega\in
   \varOmega$. Moreover, the equation\/
   $\big(\widetilde{\varphi}\ast\psi\ast\gamma\big)(t) = \langle
   f_\varphi, T^t f_\psi\rangle$ holds for arbitrary\/ $\varphi,\psi \in
   C_c (G)$ and\/ $t\in G$.
\item[\rm (b)] Let\/ $G$ have a topology with countable base. Let\/ $\oO$
   be a\/ {\rm TMDS} with ergodic probability measure\/ $m$.  Then,
   there exists a translation bounded measure\/ $\gamma$ on\/ $G$ such
   that the sequence\/ $(\frac{1}{|B_n |}\, \widetilde{ \omega^{}_{B_n}}
   \ast \omega^{}_{B_n})$ converges, in the vague topology, to\/
   $\gamma$ for\/ $m$-almost every\/ $\omega \in \Omega$, whenever\/ $(B_n)$
   is a van Hove sequence along which the Birkhoff ergodic theorem
   holds. Moreover, the equation\/ $\big(\widetilde{\varphi}\ast\psi
   \ast \gamma\big)(t) = \langle f_\varphi, T^t f_\psi\rangle$ holds
   for arbitrary\/ $\varphi,\psi \in C_c (G)$ and\/ $t\in G$.
\end{itemize}
\end{theorem}

\noindent
\noindent{\bf Remarks.}  
(a) Every LCA group with a countable base of the topology
admits a van Hove sequence along which the Birkhoff ergodic theorem
holds, as follows from recent results of Lindenstrauss \cite{Lindenstrauss},
see also Tempelman's monograph \cite{Tempelman}, in particular its
Appendix, for background material. More precisely, every
van Hove sequence is a F{\o}lner sequence, and thus contains a so-called
tempered subsequence with the desired property, compare \cite{Lindenstrauss}.
Note also that $G$ second countable implies
$\sigma$-compactness as well as metrisability of $G$.  \newline 
(b) The theorem generalises the corresponding results of
\cite{Dworkin,Martin2,Hof2}.

\smallskip

To prove Theorem~\ref{Dworkin}, we need some preparation in form of
the following results.
\begin{lemma} \label{hilf} 
   Let\/ $D$ be a dense subset of\/ $C_c (G)$. 
   Let\/ $C>0$ and a relatively compact
   open\/ $V\!$ in\/ $G$ be given.  If\/ $(\mu_\iota)$ is a net of measures
   in\/ $\MCV$ such that\/ $\mu_\iota (\varphi)$ converges for every\/
   $\varphi\in D$, then there exists a translation bounded measure\/
   $\mu \in \MM$ such that\/ $(\mu_\iota)$ converges vaguely to\/ $\mu$.
\end{lemma}
\begin{proof}  As $D$ is dense, every $\mu$ in $\MM$ is
uniquely determined by its values on $D$. Thus, all converging subnets
of $(\mu_\iota)$ have the same limit.  As $\MCV$ is compact by Theorem
~\ref{Compact}, there exist converging subnets.  Putting this
together, we arrive at the desired statement. 
\end{proof}

\smallskip
\begin{lemma} \cite[Lemma~1.2]{Martin2}  \label{Martin} 
   Let\/ $\mu,\nu$ be translation bounded measures on\/ $G$ and\/
   $(B_n)$ a van Hove sequence. Then, in the vague topology, 
   $\lim_{n\to \infty} \frac{1}{\lvert B_n \rvert}
   (\mu^{}_{B_n} \ast \nu^{}_{B_n} - \mu \ast
   \nu^{}_{B_n}) = 0$.  \hfill \qedsymbol
\end{lemma}
\begin{lemma} \cite[Lemma~1.1 (2)]{Martin2}  \label{randterm} 
   Let\/ $(B_n)$ be a van Hove sequence in\/ $G$ and\/ $\mu$ a
   translation bounded measure. Then, the sequence\/ $\big( |\mu|(B_n) /
   |B_n| \big)$ is bounded.
   \hfill \qedsymbol
\end{lemma}

\smallskip
\noindent{\it Proof of Theorem}~\ref{Dworkin}.  (a) As $\omega\in
\MCV$, Lemma~\ref{randterm} and a short calculation give a constant
$C'>0$ and and a relatively compact open $V'\subset G$ such that the
sequence $\big( ( \widetilde{\omega^{}_{B_n}} \ast 
\omega^{}_{B_n}\,) / |B_n| \big)$
is contained in $\mathcal{M}_{C',V'}$. Moreover, the set $\{\varphi
\ast \psi : \varphi,\psi \in C_c(G)\}$ is dense in $C_c (G)$ by
standard arguments involving approximate units \cite{Rudin}. Thus, by
Lemma~\ref{hilf} and Lemma~\ref{Martin}, it suffices to show
$\lim_{n\to \infty} \frac{1}{|B_n|}\big( \widetilde{\varphi}\ast \psi \ast
\widetilde{\omega^{}_{B_n}} \ast \omega^{}_{B_n}\big) (t) = \langle
f_\varphi, T^t f_\psi \rangle$ for arbitrary $\varphi,\psi\in C_c
(G)$ and $t\in G$.  By Lemma~\ref{Martin}, it suffices to show
\begin{equation*}
    \lim_{n\to \infty} \frac{1}{|B_n|} \big(\widetilde{\varphi}\ast
    \psi \ast \widetilde{\omega} \ast \omega^{}_{B_n} \big) (t) \; =
    \; \langle f_\varphi, T^t f_\psi \rangle.
\end{equation*} 
This follows by unique ergodicity and a Dworkin type calculation
\cite{Dworkin, Martin2,LMS-1}. As the details are somewhat more
involved than in the case of Delone sets, we include a sketch for the
convenience of the reader. We define $ Z_n := \big(
\widetilde{\varphi}\ast \psi \ast \widetilde{\omega} \ast
\omega^{}_{B_n} \big) (t)$. Then,
\begin{equation*}
   Z_n \; = \; \int_G \big(\widetilde{\varphi} \ast \psi\big) (t - u
   )\dd (\widetilde{\omega} \ast \omega^{}_{B_n})(u) \; = \; \int_G
   \int_G \int_G \overline{\varphi}(v -t + r)\psi(v -s) 1_{B_n} (s)\dd
   v \dd \widetilde{\omega}(r)\dd \omega(s),
\end{equation*}
where $1_{B_n}$ denotes the characteristic function of $B_n$. Using
Fubini's Theorem and sorting the terms, we arrive at
\begin{equation*} 
   Z_n \; = \; \int_G \left( \int_G \overline{\varphi} (v -t + r) \dd
  \widetilde{\omega}(r)\right) \left( \int_G \psi (v - s) 1_{B_n}
  (s)\dd \omega(s)\right)\dd v \, .
\end{equation*} 
We will now study the two terms in brackets. A short calculation shows
$ \int_G \overline{\varphi} (v -t + r)\dd \widetilde{\omega}(r) =
\overline{ f_\varphi (\alpha^{}_{t-v} \omega)}.$ As for the other
term, we consider the difference function
\begin{equation*} 
   D(v) \; := \; \int_G \psi (v - s)\, 1_{B_n} (s)\dd \omega(s) - 
                 \int_G \psi (v-s)\dd \omega(s)\, 1_{B_n} (v) \, .
\end{equation*}
Let $K$ be a compact set with $K=-K$, $0\in K$ and $\supp (\psi)
\subset K$.  Then, it is not hard to see that the difference $D(v)$
(and in fact each of its terms alone) vanishes for $v\notin B_n +
K$. Similarly, one can show that $D(v)$ is supported in $\partial^K
B_n$.  Apparently, $|D(v)|$ is bounded above by $C := 2
\|\psi\|_\infty \sup\,\{|\omega| (t +\supp(\psi) ) : t\in G\} < \infty$.
As $(B_n)$ is a van Hove sequence, we conclude that
\begin{equation*}
     0 \; \le \;
     \frac{1}{|B_n|} \int_G \left| D(v) \overline{ f_\varphi
     (\alpha^{}_{t-v} \omega)}\right |\dd v 
     \; \le \; C \|f_\varphi\|_\infty
     \frac{|\partial^K B_n|}{|B_n|} 
     \; \xrightarrow{n\to\infty} \; 0 \, .
\end{equation*}
Noting that $ \int_G \psi (v-s)\dd \omega(s) = f_\psi (\alpha^{}_{-v}
\omega)$, and putting these considerations together, we arrive at
\begin{equation*}
   \lim_{n\to \infty} \frac{ Z_n}{ |B_n| }
   \; = \; \lim_{n\to \infty}
   \frac{1}{|B_n|} \int_G \overline{ f_\varphi (\alpha^{}_{t-v} \omega)}
   f_\psi (\alpha^{}_{-v} \omega)\, 1_{B_n} (v)\dd v \, .
\end{equation*}
This yields 
\begin{eqnarray*}
   \lim_{n\to \infty} \frac{ Z_n}{|B_n|}  
   &=& \int_\varOmega \overline{ f_\varphi 
       (\alpha^{}_{t} \omega)} f_\psi (\omega)\dd m (\omega)
   \;\, = \;\, \int_\varOmega \overline{ f_\varphi (\omega)} f_\psi 
       (\alpha^{}_{-t}\omega)\dd m (\omega)\\ 
   &=& \int_\varOmega \overline{ f_\varphi (\omega)} (T^t f_\psi) 
       (\omega)\dd m (\omega)
   \;\, = \;\, \langle f_\varphi, T^t f_\psi\rangle \, , 
\end{eqnarray*}
where we used the pointwise ergodic theorem for continuous functions
on a uniquely ergodic system in the first step and $\alpha$-invariance
of $m$ in the second step. (Note that this ergodic theorem
only relies on compactness of the underlying space $\varOmega$ and
does not require separability of $G$. This can easily be seen by going through  
a  proof of this theorem as presented, e.g., in \cite[Thm.~6.19]{Wal}.) 
 
(b) This can be seen
similarly: After replacing the pointwise ergodic theorem for uniquely ergodic
systems by the Birkhoff ergodic theorem, the considerations of (a) can
be carried through to show that, for each function
$\widetilde{\varphi}\ast \psi$, there exists a set
$\varOmega_{\varphi,\psi}\subset \varOmega$ of full measure such that
\begin{equation*} 
     \frac{1}{|B_n |}\, \big(\widetilde{\varphi}\ast \psi \ast  
     \widetilde{ \omega_{B_n}} \ast \omega^{}_{B_n}\big) (0) 
     \; \xrightarrow{n\to\infty} \; \langle
     f_\varphi , f_\psi\rangle 
\end{equation*} 
for every $\omega\in\varOmega_{\varphi,\psi}$. As $G$ is second
countable, there exists a countable set $D$ in $C_c (G)$ such that $D$
and $\{\widetilde{\varphi}\ast \psi : \varphi,\psi \in D\}$ are dense
in $G$. Thus, there is a set $\varOmega_0\subset \varOmega$ of full
measure such that
\begin{equation*} 
   \frac{1}{|B_n |}\, \big(\widetilde{\varphi}\ast \psi \ast  
   \widetilde{ \omega_{B_n}} \ast \omega^{}_{B_n}\big) (0)
   \; \xrightarrow{n\to\infty} \;
   \langle f_\varphi , f_\psi\rangle 
\end{equation*} 
for all $\varphi,\psi \in D$ and all $\omega \in \varOmega_0$. By the
density of $\{\widetilde{\varphi}\ast \psi : \varphi,\psi \in D\}$ in
$C_c (G)$ and Lemma~\ref{hilf}, the vague convergence of 
$\big(\widetilde{ \omega_{B_n}} \ast \omega^{}_{B_n}\big) / |B_n|$ 
towards a translation
bounded measure $\gamma$ with $\big(\widetilde{\varphi}\ast \psi \ast
\gamma\big) (0) = \langle f_\varphi , f_\psi\rangle$ for all $\varphi, \psi
\in D$ follows.  This gives the desired vague convergence.

It remains to show the last part of the statement: As $D$ is dense in
$C_c (G)$ and $f$ is a continuous map, the formula
\begin{equation*}
  \big(\widetilde{\varphi}\ast \psi \ast \gamma\big) (0) 
  \; = \; \langle f_\varphi , f_\psi\rangle 
\end{equation*}
does not only hold for $\varphi, \psi \in D$, but for arbitrary
$\varphi, \psi \in C_c (G)$.  For $t\in G$, this implies
\begin{equation*}
   \big( \widetilde{\varphi} \ast \psi \ast \gamma \big)(t) 
   \; = \; \big( \delta_{-t} \ast \widetilde{\varphi} \ast \psi \ast
   \gamma  \big) (0) 
   \; = \; \big( \widetilde{\varphi} \ast
   (\delta_{-t} \ast \psi) \ast \gamma \big) (0) 
   \; = \;
   \langle f_{\varphi} , f_{\delta_{-t} \ast \psi}\rangle \, .
\end{equation*}
By Lemma~\ref{con}, we have 
\begin{equation*}
   f_{\delta_{-t} \ast \psi} (\omega) 
   \; = \; f_{\beta_{-t} (\psi)} (\omega) \; = \;
   f_{\psi} (\alpha^{}_{-t} \omega) 
   \; = \; T^t f_{\psi} (\omega) \, . 
\end{equation*}
Thus, we can conclude
\begin{equation*}
   \big( \widetilde{\varphi} \ast \psi \ast \gamma \big) (t) 
   \; = \; \langle f_{\varphi} , f_{\delta_{-t} \ast \psi}\rangle 
   \; = \; \langle f_{\varphi} , T^t f_\psi \rangle \, .
\end{equation*} 
This finishes the proof.  \hfill \qedsymbol

\smallskip
\noindent{\bf Remark.} By Lemma \ref{Martin}, the convergence of 
$|B_n|^{-1} \widetilde{\omega_{B_n}} \ast \omega_{B_n}$ towards 
$\gamma$ discussed in the previous  theorem  implies convergence 
of $|B_n|^{-1} \widetilde{\omega} \ast \omega_{B_n}$ towards 
$\gamma$ as well.

\smallskip
\noindent
The measure 
\begin{equation*} 
   \gamma \; = \; \gamma^{}_\omega 
   \; = \; \lim_{n\to \infty}  \frac{1}{|B_n|} 
           \, \widetilde{\omega^{}_{B_n}} \ast \omega^{}_{B_n} 
\end{equation*}
appearing in the theorem is called the \textit{autocorrelation
measure} (or \textit{autocorrelation} for short) of $\omega\in
\varOmega$. It is obviously positive definite, and hence transformable. 
By Bochner's Theorem, compare \cite[Ch.~I.4]{BF},
its Fourier transform is then a positive measure on $\widehat{G}$,
called the \textit{diffraction measure} of $\omega\in \varOmega$. We
will have to say more about autocorrelation and diffraction measures 
in the next section.

\section{Relating diffraction and dynamical spectrum}\label{Dynamical}

In this section, we show that the diffraction spectrum is equivalent
to the spectrum of a certain subrepresentation of $T$.  This type of
result is implicit in essentially every work using the so-called
Dworkin argument \cite{Dworkin,Hof,Martin2,Boris}. However, it seems
worthwhile to make this connection explicit. In fact, this is one of
the two cornerstones of our approach to the characterisation of pure
pointedness, the other being Theorem~\ref{Abstractcriterion}. A key 
ingredient in our considerations will be Proposition \ref{spectralmass} 
below.

\smallskip
We start by giving a closed formula for the autocorrelation
measure. This closed formula does not rely on any ergodicity
assumptions. Thus, via this formula, an autocorrelation measure can be
attached to any TMDS with an invariant probability measure $m$.  We
should like to mention that this is inspired by recent work of
Gou\'{e}r\'{e} \cite{Gouere-1}, who gives a closed formula in the
context of Palm measures and point processes.
\begin{prop} \label{geschlossen} 
  Let\/ $\oO$ be a\/ {\rm TMDS} with invariant probability measure\/ $m$.  
  Let a function\/ $\sigma\in C_c (G)$ be given
  with\/ $\int_G \sigma(t)\dd t = 1$.
  For\/ $\varphi\in C_c (G)$, define
\begin{equation*} 
    \gamma^{}_{\sigma,m} (\varphi) \; := \;
    \int_\varOmega \int_G f_\varphi (\alpha^{}_{-t} \overline{\omega}) 
    \sigma(t) \dd \omega(t)\dd m(\omega).
\end{equation*}
  This leads to the following assertions.
\begin{itemize}
\item[\rm (a)] The map\/ $\gamma^{}_{\sigma,m} \!:\, C_c
   (G)\longrightarrow \CC$ is continuous, i.e.,
   $\gamma^{}_{\sigma,m}\in \MM$.
\item[\rm (b)] For\/ $\varphi,\psi \in C_c (G)$, the equation\/ 
   $\big(\widetilde{\varphi}\ast \psi \ast \gamma^{}_{\sigma,m}\big) (t) 
   =  \langle f_\varphi, T^t  f_\psi \rangle$ holds. 
\item[\rm (c)] The measure\/ $\gamma^{}_{\sigma,m}$ does not depend on\/ 
   $\sigma\in C_c (G)$, provided\/ $\int_G \sigma\dd t =1$. 
\item[\rm (d)] The measure\/ $\gamma^{}_{\sigma,m}$ is positive definite. 
\end{itemize}
\end{prop}
\begin{proof}  
Note that $f_{\varphi} (\alpha_{s}
\overline{\omega})= \overline{f_{\overline{\varphi}} (\alpha_s
(\omega))}$.  \\
(a) Obviously, $| f_\varphi (\overline{\omega})| \le
\|\varphi\|_\infty \sup\,\{ |\omega|(t + \supp(\varphi) ) : t\in G\}$. As
$\sigma$ has compact support, $\gamma^{}_{\sigma,m} (\varphi)$ is then
finite. Moreover, $\varphi_\iota \longrightarrow \varphi$ implies
$\overline{\varphi_\iota} \longrightarrow \overline{\varphi}$ which, in
turn, yields $\overline{f_{\overline{\varphi_\iota}}}\longrightarrow
\overline{f_{\overline{\varphi} }}$ by continuity of $f$.  As $\sigma$
has compact support and $\varOmega\subset \MCV$, this gives
\begin{equation*} 
   \int_G \overline{f_{\overline{\varphi_\iota}} } (\alpha^{}_{-t}
   (\omega)) \sigma(t)\dd \omega(t)\; \longrightarrow \; \int_G
   \overline{f_{\overline{\varphi}} (\alpha^{}_{-t} (\omega))} \sigma(t)
  \dd \omega(t) \, , 
\end{equation*} 
uniformly on $\varOmega$. The desired continuity statement
follows. \\ 
(b) We first  show the statement for $t=0$. To this aim, we define
\begin{equation*} 
   Z \; := \; \big(\widetilde{\varphi} \ast \psi \ast 
           \gamma^{}_{\sigma,m}\big) (0) 
     \;  = \; \int_G \big(\widetilde{\varphi}\ast \psi\big) (-s) \dd 
           \gamma^{}_{\sigma,m} (s)
     \;  = \; \int_\varOmega \int_G \overline{ 
              f_{\widetilde{\widetilde{\varphi} \ast \psi  }   } }
            (\alpha^{}_{-t} \omega) \sigma(t)\dd \omega(t)\dd m (\omega).
\end{equation*}
Then, we can calculate
\begin{eqnarray*}
Z & =& \int_\varOmega \int_G \, \overline{ \int_G
\overline{(\widetilde{\varphi}\ast\psi )(s -t) }\dd \omega(s)} \,
\sigma(t)\dd \omega(t)\dd m(\omega)\\ &=& \int_\varOmega \int_G \int_G
\int_G \widetilde{\varphi} (u) \psi (s-t -u)\dd u\,\dd
\overline{\omega} (s) \sigma(t)\dd \omega(t)\dd m(\omega)\\ &=&
\int_\varOmega \int_G \int_G \int_G \widetilde{\varphi} (u + s) \psi
(-t -u)\dd u\,\dd \overline{\omega} (s) \sigma(t)\dd \omega(t)\dd
m(\omega)\\ &=& \int_G \int_\varOmega \int_G \overline{f_\varphi (
\alpha_u \omega)} \psi (-t -u) \sigma (t) \dd \omega(t)\dd
m(\omega)\dd u,
\end{eqnarray*}
where we used the translation invariance of the Haar measure in the
second last step, and Fubini's Theorem and $\int_G \widetilde{\varphi}
(u +s)\dd \overline{\omega} (s) = \overline{f_\varphi ( \alpha_u
\omega)}$ in the last step.  By the invariance of $m$, and Fubini's
Theorem together with $\int \sigma(t)\dd t = 1$, this gives:
\begin{eqnarray*}
Z &=& \int_G \int_\varOmega \int_G \overline{f_\varphi (\omega)} 
\psi (-t -u) \sigma (t)\dd (\alpha_{-u} \omega)(t)\dd m(\omega)\dd u\\
&=& \int_G \int_\varOmega \int_G \overline{f_\varphi (\omega)} \psi (-t)
\sigma (t-u )\dd  \omega(t)\dd m(\omega)\dd u\\
&=& \int_\varOmega \overline{f_\varphi(\omega)} f_\psi(\omega)\dd m(\omega)
\;\, = \;\, \langle f_{\varphi} , f_{\psi} \rangle \, .
\end{eqnarray*}
The case of arbitrary $t\in G$ can now be treated by mimicking the
last part of the proof of part (b) of Theorem \ref{Dworkin}.  \\
(c) This is
immediate from (b) and (a) as $\{ \widetilde{\varphi} \ast \psi :
\varphi,\psi \in C_c (G)\}$ is dense in $C_c (G)$.  \\(d) This is a
direct consequence of (b).
\end{proof}

\smallskip
Part (b) of the Lemma shows, in particular, that the measure
$\gamma^{}_{\sigma,m}$ equals the autocorrelation measure introduced
in the last section if $m$ is ergodic. Part (d) shows that
$\gamma^{}_{\sigma,m}$ is positive definite. Thus, by Bochner's
Theorem, see \cite{BF}, its Fourier transform is a positive measure 
on the dual group $\widehat{G}$.
\begin{definition} \label{def-auto}
  Let $\oO$ be a\/ {\rm TMDS} with invariant 
  probability measure $m$. The measure\/ $\gamma^{}_m:=
  \gamma^{}_{\sigma,m}$, where\/ $\sigma\in C_c (G)$ with\/ $\int_{G}
  \sigma\dd t =1$, is called the\/ {\em autocorrelation measure} of the
  dynamical system\/ $\oO$ with invariant measure\/ $m$. Its Fourier
  transform\/ $\widehat{\gamma^{}_m}$ is called the\/
  {\em diffraction measure}
  of the dynamical system\/ $\oO$ with invariant measure\/ $m$.
\end{definition}

We summarise the preceding considerations in the following lemma,
where we use $\psi_{\!\_}$ for the function defined by
$\psi_{\!\_} (t) = \psi(-t)$.
\begin{lemma} Let\/ $\oO$ be a\/ {\rm TMDS} with invariant probability
    measure\/ $m$. Then, there exists a unique measure on\/ $G$ assigning
    the value\/ $\langle f_\varphi, f_\psi\rangle$ to the function\/
    $\overline{\varphi} * \psi_{\!\_}$, for\/ $\varphi,\psi\in
    C_c (G)$. This measure is the autocorrelation measure\/
    $\gamma^{}_m$ of\/ $\oO$.
\end{lemma} 
\begin{proof} 
Uniqueness is clear as the set 
$\{\overline{\varphi} * \psi_{\!\_} : \varphi, \psi \in C_c (G)\}$ 
is dense in $C_c (G)$. Existence follows from 
Proposition~\ref{geschlossen}, as
\begin{equation*}
   \gamma^{}_m (\overline{\varphi} * \psi_{\!\_}) \, = \,
   \big(\widetilde{\varphi}\ast \psi \ast \gamma^{}_m\big) (0)
   \, = \, \langle f_\varphi, f_\psi\rangle.
\end{equation*}
This proves the lemma. 
\end{proof}

\noindent{\bf Remark.} 
This definition of the autocorrelation and the
diffraction of a dynamical system is to be compared with the
corresponding objects of a single measure (namely an element in
$\varOmega$) studied in the last section. In the latter case, one
faces the problem of its dependence of the measure $m$, or of the
averaging sequence $(B_n)$. It is reasonable, both mathematically and
physically, to replace this by the objects defined in 
Definition~\ref{def-auto}, at least for most aspects of the 
spectral theory connected with it.

\smallskip
Having cast the diffraction measure in an abstract context, we will
now briefly discuss the basic quantities in the spectral theory of
dynamical systems: Let $\oO$ be a TMDS. By Stone's Theorem (compare
\cite[Sec.~36D]{Loomis}), there exists a projection valued measure
\begin{equation*} 
   E_T\!:\, \mbox{Borel sets on $\widehat{G}$}\; \longrightarrow \;
   \mbox{Projections on $L^2(\varOmega,m)$} 
\end{equation*} 
with
\begin{equation*} 
   \langle f, T^t f \rangle 
   \; = \; \int_{\widehat{G}} (\widehat{s}, t)\dd 
   \langle f, E_T(\widehat{s}) f \rangle 
   \; = \; \int_{\widehat{G}} (\widehat{s},
   t)\dd  \rho^{}_f (\widehat{s}) \, , 
\end{equation*} 
where $\rho^{}_f$ is the measure on $\widehat{G}$ defined by
$\rho^{}_f (B) := \langle f, E_T (B)f\rangle$. It is then not hard to
see that $T$ has pure point spectrum (in the sense defined in
Section~\ref{Criterion}) if and only if all the measures $\rho^{}_f$,
with $f\in L^2 (\varOmega,m)$, are pure point measures.

\smallskip

To $\varphi \in C_c (G)$, we have associated the function $f_\varphi
\in L^2 (\varOmega,m)$ in the last section. It turns out that the
measure $\rho^{}_{f_\varphi}$ can be calculated in terms of the
diffraction measure. While this connection is not hard to prove, it is
underlying the main result of this section. Therefore, we isolate it in
the following proposition (compare \cite{BF,GdeL}).

\begin{prop} \label{spectralmass} 
   Let\/ $\oO$ be a\/ {\rm TMDS} with invariant probability measure\/ $m$. 
   Then, the equation\/ $\rho^{}_{f_\varphi} =
   |\widehat{\varphi}|^2 \, \widehat{\gamma_m}$ holds for every $\varphi \in
    C_c (G)$.
\end{prop}
\begin{proof}
By the very definition of $\rho_{f_\varphi}$ above, the (inverse) 
Fourier transform (on $\widehat{G}$) of  $\rho^{}_{f_\varphi}$ is\/
$t\mapsto \langle f_\varphi, T^t f_\varphi\rangle$.  
By Lemma~\ref{geschlossen}, we
have $\langle f_\varphi , T^t f_\varphi\rangle =
\big(\widetilde{\varphi} \ast \varphi \ast \gamma^{}_{m}\big) (t)$. 
Thus, taking the Fourier transform (on $G$), we infer 
$\rho^{}_{f_\varphi} = |\widehat{\varphi}|^2 \, \widehat{\gamma_m}$. 
\end{proof}

\smallskip

Note that every closed $T$-invariant subspace $\cV$ of
$L^2(\varOmega,m)$ gives rise to a representation $T|_{\cV}$ of $G$ on
$\cV$ by restricting the representation $T$ to $\cV$. The spectral
family of $T|_{\cV}$ will be denoted by $E_{T|_{\cV}}$.  With the
canonical inclusion $i^{}_{\cV} \!:\, \cV \longrightarrow
L^2(\varOmega,m)$ and projection $P^{}_{\cV}
\!:\,L^2(\varOmega,m)\longrightarrow \cV$, we obviously have
\begin{equation*}
   T|_{\cV} \; = \;  P^{}_{\cV}  \,T \, i^{}_{\cV}\quad 
   \mbox{and}\quad E_{T|_{\cV}} = P^{}_{\cV}\, E_T\,  i^{}_{\cV}.
\end{equation*}

In our setting, a  translation invariant subspace appears naturally. 
This is discussed next.
\begin{lemma} \label{AssumptionI} 
  Let\/ $\oO$ be a\/ {\rm TMDS} with invariant probability measure\/ $m$. 
  The set of functions\/ $\cU_0:= \{ f_\varphi :
  \varphi \in C_c (G)\}$ is a translation invariant subspace of\/
  $L^2(\varOmega,m)$, and so is its closure.
\end{lemma}
\begin{proof}
The first part of the statement follows from Lemma ~\ref{con} (b). The
second part of the statement is then immediate.
\end{proof}
\begin{definition} \label{definitionvonU}
  Let\/ $\cU$ be the closure of the space\/ $\cU_0$ from 
  Lemma~$\ref{AssumptionI}$ in\/ $L^2(\varOmega,m)$. 
\end{definition}

\smallskip
Before we can give a precise version of the relationship between  
$\widehat{\gamma_m}$ and $T$ we need one more definition. 
\begin{definition} 
  Let\/ $\rho$ be a measure on\/ $\widehat{G}$ and\/ $S$ be an arbitrary
  unitary representation of\/ $G$ on\/ $L^2(\varOmega,m)$.  
  Then, $\rho$ is called a {\em spectral measure}\/ 
  for\/ $S$ if the following holds for all Borel sets\/ $B$: 
  $E_S(B)=0$ if and only if\/ $\rho(B)=0$.
\end{definition}

Now, the relationship between $\widehat{\gamma_m}$ and $T=T_m$ can be
phrased as follows.
\begin{theorem}\label{Spectralmeasure} 
  Let\/ $\oO$ be a\/ {\rm TMDS} with invariant probability measure\/ $m$. 
  Then, the measure\/ $\widehat{\gamma_m}$ is a spectral measure for the
  restriction\/ $T|_{\cU}$ of\/ $T$ to\/ $\cU$.
\end{theorem}
\begin{proof}
Let $B$ be a Borel set in $\widehat{G}$. Then, obviously,
$E_{T|_{\cU}} (B)= 0$ if and only if we have $0= \langle f_\varphi,
E_{T|_{\cU}} (B) f_\varphi\rangle$ for every $\varphi \in C_c (G)$. As
$ \langle f_\varphi, E_T (B) f_\varphi\rangle = \langle f_\varphi,
E_{T|_{\cU}}(B) f_\varphi\rangle $ for every $\varphi \in C_c(G)$, by
the very definition of $\cU$ and $E_{T|_{\cU}}$, we infer that
$E_{T|_{\cU}} (B)= 0$ if and only if
\begin{equation*}  
   \langle f_\varphi, E_T (B) f_\varphi\rangle \; =\; 0 \, ,
   \quad \mbox{for every }  \varphi\in C_c (G)\, . 
\end{equation*}
By Proposition \ref{spectralmass}, we have $\rho^{}_{f_\varphi} =
|\widehat{\varphi}|^2 \widehat{\gamma_m}$ and,
in particular, 
\begin{equation*}
   \langle f_\varphi, E_T (B) f_\varphi\rangle 
   \; = \rho_{f_\varphi} (B) \; =  \int_B
     |\widehat{\varphi}|^2\dd \widehat{\gamma_m} \, .
\end{equation*} 
These considerations show that $E_{T|_{\cU}} (B)= 0$ if and only if
$0= \int_B |\widehat{\varphi}|^2\dd \widehat{\gamma_m}$ for every
function
$\varphi \in C_c (G)$.  Thus,  it remains to  be shown that 
$\widehat{\gamma_m} (B)=0$ if and only if  
$0= \int_B |\widehat{\varphi}|^2\dd \widehat{\gamma_m}$ for every 
$\varphi \in C_c (G)$. The only if part is clear. As for the converse, 
recall that the image of $L^1 (G,\dd t)$ under the Fourier transform 
separates points in $\widehat{G}$, see \cite{Rudin}. As
$C_c(G)$ is dense in in $L^1(G,\dd t)$, the same holds for the image
of $C_c (G)$ under the Fourier transform. Therefore, for every
$\widehat{s} \in \widehat{G}$, there exists a $\varphi \in C_c(G)$
with $\widehat{\varphi}(\widehat{s})\neq 0$. Thus, $\widehat{\gamma_m} 
(B\cap K)= 0$ for every compact $K$ is a direct consequence of 
$0= \int_B |\widehat{\varphi}|^2\dd \widehat{\gamma_m}$ for every 
$\varphi \in C_c (G)$. As $\widehat{\gamma_m}$ is regular, $\widehat{\gamma_m} 
(B) = 0$ follows. 
\end{proof}

We finish this section with a brief discussion of some consequences of the
above results for the definition of $\gamma_m$. We
thank the referee for useful comments on this point.

\smallskip

The following is a consequence of Proposition~\ref{spectralmass} (compare
\cite[Prop.~3]{BM} and discussion preceding it for similar considerations).
 
\begin{coro}  \label{coro1}
    Let\/ $(\varphi_\iota)$ be an approximate unit in\/ $C_c (G)$ 
    with respect to convolution, i.e., $\varphi_\iota \ast\psi
    \longrightarrow \psi$ in $C_c (G)$ for every $\psi \in C_c
    (G)$. Then, the measures $\rho_{f_{\varphi_\iota}}$ converge vaguely
    to $\widehat{\gamma_m}$.
\end{coro}
\begin{proof} 
As $(\varphi_\iota)$ is an approximate unit in $C_c (G)$, the continuous
functions $\varphi_\iota \ast \widetilde{\varphi_\iota}$, viewed as
absolutely continuous measures with respect to the Haar measure on $G$,
converge vaguely towards $\delta_0$, the unit point measure at $0\in G$. 
Thus, Levy's continuity theorem \cite[Thm 3.13]{BF} gives us compact convergence of
$|\widehat{\varphi_\iota}|^2$ towards $\widehat{\delta_0}\equiv 1$. 
This easily implies vague convergence of the measures
$|\widehat{\varphi_\iota}|^2\, \widehat{\gamma_m}$ towards
$\widehat{\gamma_m}$. As $\rho_{f_\varphi} =|\widehat{\varphi}|^2\,
\widehat{\gamma_m}$ by Proposition~\ref{spectralmass}, we infer the
statement of the corollary.
\end{proof}

Note that Corollary~\ref{coro1} gives another way to define the diffraction 
measure $\widehat{\gamma_m}$. Namely, we can define $\widehat{\gamma_m}$ to
be any accumulation point of the net $(\rho_{f_{\varphi_\iota}})$ whenever
$(\varphi_\iota)$ is an approximate unit in $C_c (G)$. The result is
unique, once $m$ is chosen.

\section{The main result}\label{Main}
In this section, we state and prove our main result. It shows
equivalence of pure point diffraction and pure point dynamical
spectrum for rather general measure theoretic dynamical systems.

\smallskip
\begin{theorem}\label{Characterization} 
  Let\/ $\oO$ be a\/ {\rm TMDS} with invariant probability measure\/ 
  $m$. Let\/ $T_m$ be the associated unitary representation of\/ $G$ 
  by translation operators and\/ $\widehat{\gamma^{}_m}$ the associated 
  diffraction measure. The following assertions are now equivalent.
\begin{itemize}
\item[\rm (a)] The diffraction measure\/ $\widehat{\gamma^{}_m}$ is
pure point.
\item[\rm (b)] The representation\/ $T_m$ has pure point spectrum. 
\end{itemize}
\end{theorem}

\smallskip
\noindent{\it Proof of Theorem~$\ref{Characterization}$.}  
(a) $\Longrightarrow$ (b). 
This is a consequence of Theorem~\ref{Abstractcriterion}. More
precisely, we will show that the vector space 
\begin{equation*}
   \cV \; := \; \{f_\varphi : \varphi \in C_c (G)\}
\end{equation*} 
satisfies assertion (b) of this theorem: 

As $\widehat{\gamma_m}$ has pure point spectrum by (a), Theorem
~\ref{Spectralmeasure} gives that $T|_{\cU}$ has pure point
spectrum, where $\cU$ is the closure of $\cV$.
Thus, in particular, $f_\varphi $ belongs to $\cHp$ for
every $\varphi \in C_c (G)$. As every element of the form $f_\varphi$
is continuous by Lemma~\ref{con}, we see that $\cV$ is indeed a
subspace of $\cHp\cap C(\varOmega)$.

It remains to be shown that $\cV$ separates points. Let $\omega^{}_1$
and $\omega^{}_2$ be two different points of $\varOmega$. Then,
$\omega^{}_1$ and $\omega^{}_2$ are different measures on
$G$. Therefore, there exists a $\varphi\in C_c (G)$ with $\omega^{}_1
(\varphi) \neq \omega^{}_2 (\varphi)$. This implies $f_{\varphi_{\!\_}}
(\omega^{}_1) \neq f_{\varphi_{\!\_}} (\omega^{}_2)$ with 
$\varphi_{\!\_}(t) := \varphi(-t)$.

\smallskip
(b) $\Longrightarrow$ (a). This is immediate from 
Theorem~\ref{Spectralmeasure}.
\hfill \qedsymbol

\section{Spectral properties determined by subrepresentations} \label{Further}

The ideas of the preceding sections can be refined to give
some further information on how spectral properties of $T$ are
determined by spectral properties of $T_{\cU}$. This concerns the
continuity of the eigenfunctions, and the set of  eigenvalues. 
While the TMDS are the application we have in mind  here, the underlying
result can be phrased rather abstractly.

\smallskip
We need a  special concept on ``density of a subspace with respect to
multiplication''.  This is defined next.
\begin{definition} A subspace\/ $\cV$ of\/ $L^2 (\varOmega,m)$ is said to
  satisfy condition\/ {\rm MD} if the set of products\/ 
  $f_1\cdot\ldots\cdot f_n$ with\/
  $n\in \NN$, $f_i \in \cV\cap L^\infty ( \varOmega,m)$ 
  or\/ $\overline{f_i} \in \cV\cap L^\infty (\varOmega, m)$, 
  $1\le i\le n$, is total in\/ $L^2 (\varOmega,m)$. 
\end{definition}
\begin{theorem}\label{Erweiterung}
Let\/ $\oO$ be a topological dynamical system over\/ $G$ with\/
$\alpha$-invariant measure\/ $m$. Let\/ $\cV$ be a closed\/ $T$-invariant
subspace of\/ $L^2 (\varOmega,m)$ satisfying\/ {\rm MD}.  If\/  $T|_\cV$ 
has pure point spectrum, then the following assertions hold:
\begin{itemize}
\item[(a)] $T$ has pure point spectrum.
\item[(b)] The group of eigenvalues of\/ $T$ is generated by the set of  
     eigenvalues of\/ $T|_\cV$. 
\item[(c)] If\/ $\cV$ has a basis consisting of
  continuous eigenfunctions of\/ $T|_\cV$, then\/ $L^2 (\varOmega,m)$ has a
  basis consisting of continuous eigenfunctions of\/ $T$, provided the
  multiplicity of each eigenvalue of\/ $T$ is at most countably infinite.
\end{itemize}
\end{theorem}

\noindent
{\bf Remarks.}  (a) The countability assumption in (c) is trivially
satisfied if the Hilbert space $L^2(\varOmega,m)$ is separable, (which
holds, e.g., if $\varOmega$ is metrisable). It is also satisfied if $\alpha$
is ergodic with respect to $m$. In this case, the multiplicity of each
eigenvalue is one. \\
(b) While we have stated the theorem for topological dynamical systems,
its proof does not use the topology on $\varOmega$. It can
therefore be carried over without changes to give the corresponding result
for measurable actions of $G$ on a measure space $\varOmega$. 
\begin{proof} (a)/(b)  Let $\cS_1$ be an orthonormal  basis of $\cV$ 
consisting of eigenfunctions of $T|_\cV$. Set 
\begin{equation*}
    \cS_2 \; := \; 
   \{ f^N : f\in \cS_1\;\:\mbox{or}\;\: 
    \overline{f}\in \cS_1, N\in \NN\},
\end{equation*}
where, for $N\in \NN$ and a function $f$, the function $f^N$ is defined 
in \eqref{cutoff}. As mentioned there,  $f^N$ is again an eigenfunction of 
$T$. However, $f^N$ need not belong to $\cV$.  Let $\cS_3$ be the set of
finite products of elements of $\cS_2$. In particular, all elements of 
$\cS_3$ are bounded functions, and the same is true of all finite 
linear combinations of elements of $\cS_3$. 

\medskip

\noindent
{\bf Claim.} 
\textit{Every finite product\/  
$f_1\cdot\ldots\cdot f_n$ with\/ $n\in \NN$, and\/ $f_i$ or\/
$\overline{f_i}$ in\/ $\cV\cap L^\infty (\varOmega,m)$, can be 
approximated arbitrarily well $($in\/ $L^2 (\varOmega,m)
\hspace{1pt})$ by finite 
linear combinations of elements of\/ $\cS_3$}.

\medskip \noindent
{\it Proof of the Claim.}  This is shown by induction. The case $n=1$ is 
simple, as $\cS_1$ is an orthonormal basis of $\cV$. Assume that the claim 
holds for fixed $n\in \NN$. As in Lemma \ref{Algebra}, we use again a variant 
of Lee, Moody and Solomyak \cite{LMS-1}. Let $\varepsilon >0$ be given. 
By the induction assumption, there exists a finite linear combination $g$ of 
elements of $\cS_3$  with 
\begin{equation*} 
   \| f_1\cdot\ldots\cdot f_n - g\|_2 
   \; \leq \; \frac{\varepsilon}{\|f_{n+1}\|}\, .
\end{equation*}
Here, $g$ is a bounded function (as all functions in $\cS_3$ are bounded). 
Thus, there exists a finite linear combination $h$ of elements in $ \cS_3$ 
with
\begin{equation*}
   \| f_{n+1} - h \|_2 \; \leq \; \frac{\varepsilon}{\|g\|_\infty} .
\end{equation*}
The proof of the claim can now be finished as in Lemma~\ref{Algebra}. 
\hfill \qed

\smallskip
The claim shows that $\cS_3$ is total in $L^2 (\varOmega,m)$, as the 
products appearing in its statement are total in $L^2 (\varOmega,m)$ 
by the density assumption MD. Now, obviously, the elements of $\cS_3$ 
are eigenfunctions of $T$ and the corresponding eigenvalues are just the 
group generated by the eigenvalues of $T|_\cV$. This proves (a) and (b). 

\smallskip
To prove (c), we consider a basis $\cS_1$ of $\cV$ consisting of 
continuous eigenfunctions of $T|_\cV$ and define 
\begin{equation*}
   \cS_4 \; := \;
  \{ f_1\cdot\ldots\cdot  f_n : n \in \NN,\; f_i \in \cS_1 
   \;\:\mbox{or}\;\: \overline{f_i} \in \cS_1\}\, .
\end{equation*}
As above, one can show that $\cS_4$ is total in $L^2 (\varOmega,m)$. 
Apparently, the elements in $\cS_4$ are continuous  eigenfunctions of $T$. 
Moreover, by general principles, eigenfunctions belonging to different 
eigenvalues are orthogonal. We now apply the Gram-Schmidt orthogonalisation 
procedure in each eigenspace, compare \cite[Sec.~3.1.13]{Ped}.
This is possible because the multiplicity of each 
eigenspace is at most countably infinite by (c). As a result, we obtain a 
basis of eigenfunctions which are continuous. (Note that Gram-Schmidt deals 
only with finite sums in each step and therefore  does not destroy 
continuity. ) 
\end{proof}

The preceding considerations can be applied to any TMDS. This will 
briefly be discussed next. To apply Theorem \ref{Erweiterung}, we need the 
following reformulation of previous results.  
\begin{prop} \label{prop7}
  Let\/ $\oO$ be a\/ {\rm TMDS} with invariant probability
  measure\/ $m$, let\/ $\cU$  be the space introduced in
  Definition~$\ref{definitionvonU}$ and denote the
  characteristic function of\/ $\varOmega$ by\/ $1_\varOmega$. 
  Then, the subspace\/
  $\cV := \cU + \{c \, 1_\varOmega : c\in \CC\}$ 
  is closed, invariant and satisfies assumption\/ {\rm MD}.
\end{prop}
\begin{proof} 
  The subspace $\cS := \{c\, 1_\varOmega  : c\in \CC\}$ is one-dimensional.
  So, as $\cU$ is closed, $\cV = \cU +
  \cS$ is closed as well. As $\cU$ and $\cS$ are $T$-invariant, so is
  $\cV$. It remains to be shown that MD is satisfied.  This is a
  consequence of the Stone-Weierstra{\ss} Theorem as $\{f_\varphi : \varphi 
  \in C_c (G)\}$ separates points (see the proof of 
  Theorem~\ref{Characterization}). 
\end{proof}

It is possible to base the proof of our main result, 
Theorem~\ref{Characterization}, on Theorem~\ref{Erweiterung} and 
Proposition~\ref{prop7}.  Here, our focus is in a somewhat different direction. 
\begin{theorem}
 Let\/ $\oO$ be a\/ {\rm TMDS} with invariant probability measure\/ 
  $m$, $T_m$ be the corresponding unitary representation of\/ $G$ 
  by translation operators, and\/ $\widehat{\gamma^{}_m}$ the associated 
  diffraction measure. Let $\cU$ be the space of functions defined in
  Definition~$\ref{definitionvonU}$.
  If\/ $\widehat{\gamma^{}_m}$ is a pure point measure, the following 
  assertions hold.
\begin{itemize}
\item[\rm (a)] The group of eigenvalues of\/ $T_m$ is generated by the set
    of points in\/ $\widehat{G}$ of positive $\widehat{\gamma^{}_m}$ measure, 
    i.e., the points\/
    $\hat{s}$ with\/ $\widehat{\gamma^{}_m}(\{\hat{s}\}) > 0$. 
\item[\rm (b)] If\/ $\cU$ has a basis of continuous eigenfunctions of\/ $T_m$, 
    then so has\/ $L^2 (\varOmega,m)$, provided the multiplicity of each 
    eigenvalue is at most countably infinite. 
\end{itemize}
\end{theorem}
\begin{proof} As $\widehat{\gamma^{}_m}$ is a pure point measure, 
$T|_{\cU}$ has pure point spectrum by Theorem~\ref{Spectralmeasure}.
Set $\cV:= \cU + \{c\, 1_\varOmega  : c\in \CC\}$. As $1_\varOmega$ is
obviously a (continuous) eigenfunction of $T$ (to the eigenvalue $1$),
$T_{\cV}$  has pure point spectrum as well.  Moreover, by 
Proposition~\ref{prop7}, $\cV$  is invariant, closed and satisfies 
assumption MD.  Thus, the conditions of Theorem \ref{Erweiterung} are 
satisfied, and our assertions follow. 
\end{proof}

\medskip

\subsection*{Acknowledgements}
It is our pleasure to thank Robert V.\ Moody and Peter Stollmann
for their cooperation. DL would also like to thank Jean-Baptiste
Gou\'{e}r\'{e} for very stimulating discussions.  This work was
supported by the German Research Council (DFG). We are grateful to the
Erwin Schr\"odinger International Institute for Mathematical Physics
in Vienna for support during a stay in winter 2002/2003, where this
manuscript was completed.

\medskip

\begin{appendix}

\section{The local rubber topology is  a Fell topology}
The aim of this appendix is to show that the topology introduced in
Section \ref{Delone} is a special case of a topology introduced by Fell in
\cite{Fell} on the closed subsets of an arbitrary locally compact
space.

\medskip

We start by recalling the definition of Fell's topology: The locally
compact space in question is $G$. For a compact set $C$ in $G$ and a
finite family $\CalF$ of open sets in $G$, we define 
$\UFell$ by
\[
   \UFell : \{C\in \cC : \varLambda \cap C = \varnothing \;\:\mbox{and}\;\:  
   \varLambda \cap A\neq \varnothing \:\;\mbox{for every $A\in \CalF$} \}.
\]
The family of all $\UFell$ with $C$ compact in
$G$ and $\CalF$ a finite family of open sets in $G$ is a basis of the
Fell topology.   This is a typical example of a so-called ``hit and
miss'' topology, where $\UFell$ consists of all closed sets which hit
the sets of $\CalF$ and miss the set $C$.

\smallskip

This topology agrees with the one introduced in Section \ref{Delone}, as
follows from the next lemma.

\begin{lemma} 
   {\rm (a)} Let\/ $C \subset G$ compact and\/ $\CalF$ a finite family of 
   open subsets of\/ $G$ be given with\/ $\UFell \neq \varnothing$. Then, 
   there exists a closed set $H$ in $G$, a  compact\/ $K\subset G$  and an 
   open neighbourhood\/ $V$ of\/ $0\in G$ with 
\[
     \Uwir \, \subset \; \UFell.
\] 
   {\rm (b)} Let a closed set $H$ in $G$, a compact subset\/ $K\subset G$ and an 
   open neighbourhood\/ $V$ of\/ $0\in G$ be given. Then, there exists a
   compact\/ $C \subset G$ and a finite family of open sets\/ $\CalF$ in\/ $G$ 
   with\/ $\UFell \neq \varnothing$ and 
\[
     \UFell \, \subset \; \Uwir.
\]
\end{lemma}
 \begin{proof} 
(a)  By $\UFell\neq\varnothing$, we have  $A\setminus C\neq
 \varnothing$ for every  $A\in \CalF$.  Thus, in  every  $A\in
 \CalF$ there exists an  $x_A\in A\setminus C$. As  $A\setminus C$ is open 
and $\CalF$  is a finite family, we can find a neighbourhood  $V$ of $0$ in $G$
 with $V= -V$ such that 
\begin{equation}\label{defxa}
    x_A + \overline{V} \, \subset \, A\setminus C
    \; \quad\mbox{for every  $A\in\CalF$}.
\end{equation}
Define $H:=\{x_A : A\in \CalF\} \quad \mbox{and}\quad K:=C\cup (H + \overline{V})$.
Then, $K$ is the disjoint union of $C$ and $H + \overline{V}$ by the very 
construction of $H$. 

\smallskip
Now, let an arbitrary   $L\in \Uwir$ be given. We have to show that $L\in \UFell$: \\
By $L\in \Uwir$, we have  $L\cap K \subset H + V \subset H + \overline{V}
\subset K\setminus C$ and, as $C\subset K$,  this implies  
\[
    L\cap C \; = \;  L \cap C \cap K \; = \; \varnothing.
\]
Moreover, for every $A\in \CalF$, we have $x_A\in H =
H \cap K \subset L + V$ and therefore  $ (x_A - V) \cap L \neq
\varnothing$. By  $V= -V$ and \eqref{defxa}, this implies 
\[
   \varnothing \; \neq \;
   (x_A + V) \cap L \; \subset \;  (A\setminus C)  \cap L 
   \; \subset \; A\cap L
\]
for every  $A\in \CalF$.  These considerations show $L\in \UFell$.  
As $L\in \Uwir$ was arbitrary, we infer $\Uwir\subset \UFell$. 

\medskip

(b) Let  $W$ be an open neighbourhood of $0$ in $G$ with  $W = - W $ and $W + W
\subset V$. Define $C := K\setminus (H + W)$,
where $H$ and $K$ are given by assumption.
As $K\cap H$ is compact, there exist $t_1, \ldots, t_n\in G$ with 
\begin{equation}\label{comp} 
    K\cap H  \; \subset \; \bigcup_{i=1}^n (t_i + W) 
    \quad \mbox{and} \quad 
    (t_j + W)\cap (K\cap H) \, \neq \, \varnothing
     \quad \mbox{for $1 \le j\le n$.}
\end{equation}
Set $\CalF := \{ t_j + W : 1\le j \le n\}$.

\smallskip
Let now $L\in \UFell$ be arbitrary. We have to show that $L\in \Uwir$: \\
By  $L\in \UFell$ and the  definition of $C$, we have
$ \varnothing = L\cap C = L \cap (K\setminus (H + W))$, so
\[
   L \cap K = L\cap (K \cap (H+ W)) \cup L\cap (K \setminus (H + W))
  = L\cap (K \cap (H + W)) \subset H + W\subset H + V.
\]
By $L \cap (t_j + W) \neq \varnothing$, $1\le j\le n$,
and $W - W \subset V$, we also have $t_j + W \subset
L + V$. Combined with \eqref{comp}, this implies 
\[
    H \cap K \; \subset \; \bigcup_{j=1}^n (t_j + W) 
    \; \subset \; L + V
\]
and we infer $L\in \Uwir$. As $L\in \UFell$ was arbitrary, the inclusion 
$\UFell\subset \Uwir$ is established. Moreover, $\UFell$ is not empty, 
as it obviously contains $H$. 
\end{proof}

\end{appendix}


\clearpage
\begin{thebibliography}{99}

\bibitem{MB}
M.~Baake,
\textit{A guide to mathematical quasicrystals},
in: {\em Quasicrystals -- An Introduction to Structure,
Physical Properties and Applications}, eds.\ J.-B.\ Suck, P.\
H\"aussler and M.\ Schreiber, Springer, Berlin (2002),
pp.\ 17--48; math-ph/9901014.

\bibitem{BG} 
M.~Baake and U.~Grimm, 
\textit{A guide to quasicrystal literature}, in: 
\textit{Directions in Mathematical Quasicrystals}, 
eds.\ M.\ Baake and R.\ V.\ Moody, CRM Monograph
Series, vol.\ 13, AMS, Providence, RI (2000), pp.\ 371--373.

\bibitem{BH}
M.~Baake and M.~H\"offe,
\textit{Diffraction of random tilings:\ some rigorous results},
J.\ Stat.\ Phys.\ {\bf 99} (2000) 219--261; math-ph/9904005.

\bibitem{BL}
M.~Baake and D.~Lenz, 
\textit{Deformation of Delone dynamical systems and topological conjugacy},  
preprint (2004).

\bibitem{BM}
M.~Baake and R.~V.~Moody,
\textit{Weighted Dirac combs with pure point diffraction}, 
J.\ Reine Angew.\ Math.\ (Crelle), in press;
math.MG/0203030.

\bibitem{BMP}
M.~Baake, R.~V.~Moody and P.~A.~B.~Pleasants, 
\textit{Diffraction from visible lattice points and $k$-th power free 
integers}, Discr.\ Math.\ {\bf 221} (2000) 3--42;
math.MG/9906132. 

\bibitem{BF}
C.~Berg and G.~Forst,
\textit{Potential Theory on Locally Compact Abelian Groups},
Springer, Berlin (1975).

\bibitem{BT}
E.~Bombieri and J.~E.~Taylor, 
\textit{Which distributions of matter diffract? An initial 
investigation}, J.\ Physique Colloque C-3 {\bf 47} (1986) C3-19 -- C3-28.

\bibitem{Duneau} 
G.~Bernuau and M.~Duneau, 
\textit{Fourier analysis of deformed model sets}, in: 
\textit{Directions in Mathematical Quasicrystals}, 
eds.\ M.\ Baake and R.\ V.\ Moody, CRM Monograph
Series, vol.\ 13, AMS, Providence, RI (2000), pp.\ 43--60.

\bibitem{Cowley}
J.~M.~Cowley,
\textit{Diffraction Physics}, 3rd ed., 
North-Holland, Amsterdam (1995). 

\bibitem{Dworkin}
S.~Dworkin, 
\textit{Spectral theory and $X$-ray diffraction},
J.\ Math.\ Phys.\ {\bf 34} (1993) 2965--2967.

\bibitem{EM}
A.~C.~D. van Enter and J.~Mi\c{e}kisz, 
\textit{How should one define a (weak) crystal?}, 
J.\ Stat.\ Phys.\ {\bf 66} (1992) 1147--1153.

\bibitem{Fell} 
J.~M.~G.~Fell, 
\textit{A Hausdorff topology for the
closed subsets of a locally compact non-Hausdorff space}, 
Proc.\ Amer.\ Math.\ Soc.\ {\bf 13} (1962) 472--476.

\bibitem{GdeL}
J.~Gil de Lamadrid and L.~N.~Argabright, 
\textit{Almost Periodic Measures},
Memoirs of the AMS, vol.\ 428, AMS, Providence, RI (1990).

\bibitem{Gouere-1}
J.-B.~Gou\'{e}r\'{e}, 
\textit{Diffraction and Palm measure of point processes}, 
talk given at the meeting \textit{Open problems in quasicrystals}, 
at CIRM, Luminy, October 2002;
to appear in: C.\ R.\ Acad.\ Sci.\ Paris; 
preprint math.PR/0208064.  

\bibitem{Gouere-3}
J.-B.~Gou\'{e}r\'{e}, 
\textit{Quasicrystals and almost periodicity}, 
preprint math-ph/0212012. 

\bibitem{Hof}
A.~Hof,
\textit{On diffraction by aperiodic structures},
Commun.\ Math.\ Phys.\ {\bf 169} (1995) 25--43.

\bibitem{Hof2}
A.~Hof,
\textit{Uniform distribution and the projection method},
in: \textit{Quasicrystals and Discrete Geometry},
ed.\ J.\ Patera, Fields Institute Monographs, vol.\ 10,
AMS, Providence, RI (1998), pp.\ 201--206.

\bibitem{Ni}
T.~Ishimasa, H.~U.~Nissen and Y.~Fukano,
\textit{New ordered state between crystalline and amorphous
in\/ {\tt Ni-Cr}\/ particles},
Phys.\ Rev.\ Lett.\ {\bf 55} (1985) 511--513.

\bibitem{Kel}
J.~L.~Kelley,
\textit{General Topology},
van Nostrand, Princeton, NJ (1955); 
reprint, Springer, New York (1975).

\bibitem{KN}
P.~Kramer and R.~Neri, 
\textit{On periodic and nonperiodic space fillings of $E^n$ obtained 
by projection}, Acta Cryst.\ A {\bf 40}
(1984) 580--587;  Erratum: Acta Cryst.\ A {\bf 41} (1985) 619.

\bibitem{LMS-1}
J.-Y.~Lee, R.~V.~Moody and B.~Solomyak,
\textit{Pure point dynamical and diffraction spectra},
Annales H.\ Poincar\'{e} {\bf 3} (2002) 1003--1018;
mp\_arc/02-39.

\bibitem{LS}
D.~Lenz and P.~Stollmann, 
\textit{Delone dynamical systems and associated random operators},
to appear in: OAMP Proceedings, Constanta 2001;
preprint math-ph/0202042.

\bibitem{Lindenstrauss}
E.~Lindenstrauss, 
\textit{Pointwise theorems for amenable groups}, 
Invent.\ Math.\ {\bf 146} (2001) 259--295.

\bibitem{Loomis}
L.~H.~Loomis, 
\textit{An Introduction to Abstract Harmonic Analysis},  
van Nostrand, Princeton. NJ (1953).

\bibitem{Moody}
R.~V.~Moody,
\textit{Model sets:~A survey}, in:\ 
{\em From Quasicrystals to More Complex Systems}, 
eds.\ F.\ Axel, F.\ D\'enoyer and J.\ P.\ Gazeau,
Springer, Berlin (2000), pp.\ 145--166;
math.MG/0002020.

\bibitem{Moody2002}
R.~V.~Moody, private communication (2002). 

\bibitem{Ped}
G.~K.~Pedersen, 
\textit{Analysis Now}, 
Springer, New York (1989); Revised \ printing (1995).

\bibitem{Pet}
K.~E.~Petersen,
\textit{Ergodic Theory},
Cambridge University Press, Cambridge (1989).

\clearpage
\bibitem{RS}
M.~Reed and B.~Simon,
\textit{Methods of Modern Mathematical Physics. I: Functional Analysis},
2nd ed., Academic Press, San Diego, CA (1980).

\bibitem{Reiter}
H.~Reiter and J.~D.~Stegeman,
\textit{Classical Harmonic Analysis and Locally Compact Groups},
Clarendon Press, Oxford (2000).

\bibitem{Rudin}
W.~Rudin,
\textit{Fourier Analysis on Groups}, 
Wiley, New York (1962); reprint (1990).

\bibitem{Martin1}
M.~Schlottmann,
\textit{Cut-and-project sets in locally compact Abelian groups}, in:
\textit{Quasicrystals and Discrete Geometry}, ed.\ J.\ Patera,
Fields Institute Monographs, vol.\ 10,
AMS, Providence, RI (1998), pp.\ 247--264.

\bibitem{Martin2}
M.~Schlottmann,
\textit{Generalized model sets and dynamical systems}, in:
\textit{Directions in Mathematical Quasicrystals},
eds.\ M.\ Baake and R.\ V.\ Moody, CRM Monograph Series,
vol.\ 13, AMS, Providence, RI (2000), pp.\ 143--159.

\bibitem{SBGC} 
D.~Shechtman, I.~Blech, D.~Gratias and J.~W.~Cahn,
\textit{Metallic phase with long-range orientational order and no
translational symmetry}, 
Phys.\ Rev.\ Lett.\ {\bf 53} (1984) 1951--1953.

\bibitem{Boris}
B.~Solomyak, 
\textit{Spectrum of dynamical systems arising from Delone sets},
in: \textit{Quasicrystals and Discrete Geometry}, ed.\ J.\ Patera,
Fields Institute Monographs, vol.\ 10, 
AMS, Providence, RI (1998), pp.\ 265--275.

\bibitem{Boris2}
B.~Solomyak,
\textit{Dynamics of self-similar tilings},
Ergod.\ Th.\ \& Dynam.\ Syst.\ {\bf 17} (1997) 695--738;
Erratum: Ergod.\ Th.\ \& Dynam.\ Syst.\ {\bf 19} (1999) 1685. 

\bibitem{Tempelman}
A.~Tempelman,
\textit{Ergodic Theorems for Group Actions},
Kluwer, Dordrecht (1992).

\bibitem{Wal}
P.~Walters,
\textit{An Introduction to Ergodic Theory},
Springer, New York (1982).

\end{thebibliography}
\end{document}